\documentclass[12pt,reqno]{amsart}
\usepackage{amssymb,delarray}
\usepackage{amsfonts}
\usepackage{epsfig}
\usepackage[all]{xy}
\def\leq{\leqslant}
\def\geq{\geqslant}

\newtheorem{thm}{Theorem}[section]
\newtheorem{lem}%[thm]
{Lemma}[section]
\newtheorem{prop}%[thm]
{Proposition}[section]
\newtheorem{claim}%[thm]
{Claim}[section]

\newtheorem{cor}%[thm]
{Corollary}[section]
\newtheorem{rem}%[thm]
{Remark}[section]

{\catcode`\@=11
\gdef\n@te#1#2{\leavevmode\vadjust{%
 {\setbox\z@\hbox to\z@{\strut#1}%
  \setbox\z@\hbox{\raise\dp\strutbox\box\z@}\ht\z@=\z@\dp\z@=\z@%
  #2\box\z@}}}
\gdef\leftnote#1{\n@te{\hss#1\quad}{}}
\gdef\rightnote#1{\n@te{\quad\kern-\leftskip#1\hss}{\moveright\hsize}}
\gdef\?{\FN@\qumark}
\gdef\qumark{\ifx\next"\DN@"##1"{\leftnote{\rm##1}}\else
 \DN@{\leftnote{\rm??}}\fi{\rm??}\next@}}

\begin{document}
\baselineskip=13.7pt plus 2pt %13.7pt plus 2pt

\title[Factorization semigroups] {Factorization semigroups and irreducible components of
Hurwitz space}
\author[Vik.S. Kulikov]{Vik.S. Kulikov}

\address{Steklov Mathematical Institute}
 \email{kulikov@mi.ras.ru}
%\curraddr{}
\dedicatory{} \subjclass{}
\thanks{This
research was partially supported by grants of NSh-9969.2006.1 and
RFBR 08-01-00095. }
%\date{This version: February 2004. }
\begin{abstract}
We introduce a natural structure of a semigroup (isomorphic to a
factorization semigroup of the unity in the symmetric group) on the
set of irreducible components of Hurwitz space of marked degree $d$
coverings of $\mathbb P^1$ of fixed ramification types. It is proved
that this semigroup is finitely presented. The problem when
collections of ramification types define uniquely the corresponding
irreducible components of the Hurwitz space is investigated. In
particular, the set of irreducible components of the Hurwitz space
of three-sheeted coverings of the projective line is completely
described.
\end{abstract}

\maketitle
%\pagenumbering{roman}
\setcounter{tocdepth}{1}
%\setcounter{tocdepth}{1}
%\tableofcontentsSymplectic isotopy of sections

%\setcounter{section}{6}

\def\st{{\sf st}}

%\tableofcontents

%\setcounter{section}{-1}

\section*{Introduction} \label{introduc}
Usually, to investigate the Hurwitz space $\text{HUR}_d(\mathbb
P^1)$ of degree $d$ coverings of the projective line $\mathbb
P^1:=\mathbb C\mathbb P^1$, the following approach is used. A Galois
group $G$ of the coverings, the number $b$ of branch points, and the
types of local monodromies (that is, collections consisting of $b$
conjugacy classes of $G$) are fixed, and after that the set of
collections of representatives of these conjugacy classes is
investigated  up to, so called, Hurwitz moves (see, for example,
\cite{Cl} -- \cite{Ka}). There are several problems (for example, to
describe the set of plane algebraic curves up to equisingular
deformation or, more generally, to describe the set plane
pseudoholomorphic curves up to symplectic isotopy, to describe the
set of symplectic Lefschetz pencils up to diffeomorphisms, and so
on) in which also resembling objects naturally arise, namely, finite
collections of elements of some group considering up to Hurwitz
moves (see, for example, \cite{M-T} -- \cite{A}). (In the case of
plane algebraic and pseudoholomorphic curves, to obtain such
collections, one should choose a pencil of (pseudo)lines to obtain a
fibration over $\mathbb P^1$.) As it was shown in  \cite{Ku}, there
is natural structure of semigroups on the sets of such collections
considered up to Hurwitz moves, namely, so called, factorization
semigroups over groups. Moreover, if we consider such fibrations not
only over the hole $\mathbb P^1$ but also over the disc $D_R=\{ z\in
\mathbb C\mid |z|\leq R\}$, then this semigroup structure has a
natural geometric meaning (see \cite{Ku}).

In section 1 of this article, we give basic definitions and
investigate properties of factorization semigroups over finite
groups. In particular, we prove that the factorization semigroups of
the unity in finite groups are finitely presented, and also we
investigate the problem when an element of factorization semigroup
is defined uniquely by its type and product.

In section 2, factorization semigroups over symmetric groups
$\mathcal S_d$ are considered more closely. Here we prove a
stabilization theorem and completely describe the factorization
semigroup of the unity in $\mathcal S_3$.

In section 3, we introduce a natural structure of a semigroup (a
factorization semigroup of the unity in symmetric group) on the set
of irreducible components of Hurwitz space of marked  degree $d$
coverings of $\mathbb P^1$ with fixed ramification types and we show
that this structure induces a semigroup structure on the set of
irreducible components of the Hurwitz space $\text{HUR}_d^{G}$ of
Galois coverings of $\mathbb P^1$ with Galois group $G$ having no
outer automorphisms. Also, the results, obtained in sections 1 and
2, are applied to the problem when the irreducible components of the
$\text{HUR}_d(\mathbb P^1)$ are defined uniquely by collections of
types of local monodromies of the coverings.
\newline {\bf Acknowledgement.} Part of this work was done at MPIM,
Bonn. I would like to thank this institution for hospitality.

\section{Semigroups over groups}

\subsection{Factorization semigroups} \label{semigr} A collection
$(S,G,\alpha,\lambda) $, where $S$ is a semigroup, $G$ is a group,
and $\alpha :S\to G$, $\lambda :G\to \text{Aut}(S)$ are
homomorphisms, is called {\it a semigroup $S$ over a group $G$} if
for all $s_1,s_2\in S$ we have
$$s_1\cdot s_2=\rho
(\alpha (s_1))(s_2)\cdot s_1= s_2\cdot\lambda (\alpha (s_2))(s_1),$$
where $\rho (g)=\lambda(g^{-1})$.

Let $(S_1,G_1,\alpha_1,\lambda_1)$ and $(S_2,G_2,$
$\alpha_2,\lambda_2)$ be two semigroups over, respectively, groups
$G_1$ and $G_2$. We call a pair $(h_1,h_2)$ of homomorphisms
$h_1:S_1\to S_2$ and $h_2: G_1\to G_2$ {\it a homomorphism of
semigroups over groups} if
\begin{itemize}
\item[{($i$)}] $h_2\circ \alpha_{1}=\alpha_{2}\circ h_1$,
\item[($ii$)] $\lambda_{2}(h_2(g))(h_1(s))=h_1(\lambda_{1}(g))(s)$
for all $s\in S_1$ and all $g\in G_1$.
\end{itemize}

The {\it factorization semigroups} defined below constitute the
principal, for our purpose, examples of semigroups over groups.

Let $O\subset G$ be a subset of a group $G$ invariant under the
inner automorphisms.  We call  the pair $(G,O)$ an {\it equipped
group}. Let us associate to the set $O$ an alphabet $X=X_O=\{ x_g
\mid g\in O\}$ and for each pair of letters $x_{g_1}, x_{g_2}\in X$,
$g_1\neq g_2$ denote by $R_{g_1,g_2;l}$ and $R_{g_1,g_2;r}$ the
following relations: $R_{g_1,g_2;l}$ has the form
\begin{equation} \label{rel1} x_{g_1}\cdot
x_{g_2}=x_{g_2}\cdot\,x_{g_2^{-1}g_1g_2}\end{equation} if
$g_2\ne\bold 1$ and $x_{g_1}\cdot x_{\bold 1}=x_{g_1}$ if $g_2=\bold
1$, and $R_{g_1,g_2;r}$ has the form
\begin{equation} \label{rel2}
x_{g_1}\cdot x_{g_2}=x_{g_1g_2g_1^{-1}}\cdot x_{g_1}\end{equation}
if $g_1\ne\bold 1$ and $x_{\bold 1}\cdot x_{g_2}=x_{g_2}$ if
$g_1=\bold 1$.

Put
$$
\mathcal{R}=\{ R_{g_1,g_2;r},
 R_{g_1,g_2;l} \mid (g_1,g_2)\in O\times O,\, g_1\ne g_2 \} ,
$$
and, with the help of the set of relations $\mathcal R$, define a
semigroup
$$S(G,O)=\langle\, x_g\in X \mid R\in \mathcal{R}\, \rangle $$
which is called the {\it factorization semigroup} of $G$ with
factors in $O$.

Introduce also a homomorphism $\alpha :S(G,O)\to G$ given by $\alpha
(x_g)=g$ for each $x_g\in X$ and call it the {\it product
homomorphism}.

Next, we define an action $\lambda$ of the group $G$ on the set $X$
as follows:
$$ x_a\in X\mapsto \lambda (g)(x_a)=x_{g^{-1}ag}\in X.$$
As is easy to see, the above relation set $\mathcal{R}$ is preserved
by the action $\lambda$. Therefore  $\lambda$ defines a homomorphism
$\lambda : B\to \text{Aut} (S(G,O))$ (the {\it conjugation action}).
The action $ \lambda (g)$ on $S(G,O)$ is called the {\it
simultaneous conjugation } by $g$. Put $\lambda _{S} =\lambda
\circ\alpha$ and $\rho _{S} =\rho \circ\alpha$.

\begin{claim} {\rm (\cite{Kh-Ku})} \label{cl1}
For all $s_1,\, s_2\in S(G,O)$ we have
\[
s_1\cdot s_2=s_2\cdot \lambda_S(s_2)(s_1)=\rho_S(s_1)(s_2)\cdot s_1.
\]
\end{claim}

It follows from Claim \ref{cl1} that $(S(G,O),G, \alpha, \lambda)$
is a semigroup over $G$. When $G$ is fixed, we abbreviate $S(G,O)$
to $S_O$. By $x_{g_1}\cdot .\, .\, .\, \cdot x_{g_n}$ we denote the
element in $S_O$ defined by a word $x_{g_1}\dots x_{g_n}$.

Notice that $S: (G,O)\mapsto (S(G,O),G, \alpha, \lambda)$ is a
functor from the category of the equipped groups to the category of
the semigroups over groups. In particular, if $O_1\subset O_2$ are
two  sets invariant under the inner automorphisms of $G$, then the
identity map $id: G\to G$ defines an embedding
$id_{O_1,O_2}:S(G,O_1)\to S(G,O_2)$. So that, for each group $G$,
the semigroup $S_G=S(G,G)$ is an {\it universal factorization
semigroup} of elements in $G$, which means that each semigroup $S_O$
over $G$ is canonically embedded in $S_G$ by $id_{O,G}$.

Let $\Gamma$ be a subgroup of $G$. Denote by $S_O^{\Gamma}= \{ s\in
S_O\mid \alpha (s)\in \Gamma\}$. Obviously, $S_{O}^{\Gamma}$ is a
subsemigroup of $S_O$ and it coincides with the image of semigroup
$S(\Gamma,O\cap \Gamma)$ under the homomorphism induced by the
inclusion $\Gamma \hookrightarrow G$. In particular, if $G_{O}$ is
the subgroup of $G$ generated by the elements of the image of
$\alpha :S_O\to G$, then $S(G_O,O)\simeq S_O^{G_O}$.

If $\Gamma=\{ {\bf{1}}\}$, then the semigroup $S_O^{\bf{1}}$ will be
denoted by $S_{O,{\bf{1}}}$ and for each subgroup $\Gamma$ of $G$ we
denote $S_{O,{\bf{1}}}^{\Gamma}=S_{O,{\bf{1}}}\cap S_{O}^{\Gamma}$.

A group $G$ acts on itself by inner automorphisms, that is, for any
group $G$ there is a natural homomorphism $h:G\to Aut(G)$ (the
action of the image $h(g)=a$ of an element $g$ on $G$ is given by
$(g_1)a=g^{-1}g_1g$ for all $g_1\in G$).  It is easy to see that the
homomorphism $h$ defines  on $S_G$ a structure of a semigroup over
$A=Aut(G)$, where the homomorphism $\alpha_A:S_G\to Aut(G)$ is the
composition $h\circ \alpha$ and an element $a\in Aut(G)$ acts on
$S_G$ by the rule $x_g\mapsto x_{(g)a}$. It is easy to see that the
subsemigroup $S_{G,\bf{1}}$ is invariant under the action of
$Aut(G)$ on $S_G$. Therefore $S_{G,\bf{1}}$ also can be considered
as a semigroup over $Aut(G)$.

To each element $s=x_{g_1}\cdot .\, .\, .\, \cdot x_{g_n}\in S_O$,
$g_i\neq \bold 1$, let us associate a number $ln(s)=n$ called the
{\it length} of $s$. It is easy to see that $ln: S_O\to \mathbb
Z_{\geq 0}=\{ {\bf{a}}\in \mathbb Z\mid {\bf{a}}\geq 0\}$ is a
homomorphism of semigroups.

For each $s=x_{g_1}\cdot\, .\, . \, . \, \cdot x_{g_n}\in S_O$
denote by $G_s$ the subgroup of $G$ generated by the images $\alpha
(x_{g_1})=g_1,\dots ,\alpha (x_{g_n})=g_n$ of the factors $x_{g_1},
\dots , x_{g_n}$.

\begin{claim} \label{cc} The subgroup $G_s$ of $G$ is well defined, that is, it does not depend on a presentation of
$s$ as a product of  generators $x_{g_i}\in X_O$.
\end{claim}

The proof of Claim \ref{cc} and the following proposition is very
simple and therefore it will be omitted.

\begin{prop} {\rm (\cite{Kh-Ku})} \label{simple} Let $(G,O)$ be an equipped group and let $s\in S_O$. We have
\begin{itemize} \item[($1$)] $\text{ker}\, \lambda $ coincides with the centralizer
$C_{O}$ of the group $G_{O}$ in $G$;
\item[($2$)] if
$\alpha (s)$ belongs to the center $Z(G_s)$ of $G_s$, then for each
$g\in G_s$ the action $\lambda (g)$ leaves fixed the element $s\in
S_O$;
\item[($3$)] if $\alpha (s\cdot x_g)$ belongs to the center $Z(G_{s\cdot x_g})$ of
$G_{s\cdot x_g}$, then $s\cdot x_g=x_g\cdot s$,
\item[($4$)] if $\alpha (s)=\bf{1}$, then $s\cdot s'=s'\cdot s$
for any $s'\in S_G$.
\end{itemize}
\end{prop}

\begin{claim}\label{commut} For any equipped group $(G,O)$ the
semigroup $S_{O,\bf{1}}$ is contained in the center of the semigroup
$S_G$ and, in particular, it is a commutative subsemigroup.
\end{claim}
\proof It follows from Proposition \ref{simple} ($4$). \qed \\

It is easy to see that if $g\in O$ is an element of order $n$, then
$x_g^n\in S_{O,\bf{1}}$.
\begin{lem} \label{fried}  Let $s\in
S_{O,\bf{1}}$ and $s_1\in S_O$ be such that $G_{s_1}=G_O$. Then
\begin{equation} \label{oo} s\cdot s_1= \lambda(g)(s)\cdot s_1
\end{equation}
for all $g\in G_O$.

In particular, if $C\subset O$ is  a conjugacy class of elements of
order $n_C$ and $s\in S_O$ is such that $G_s=G$, then for any
$g_1,g_2\in C$ we have \begin{equation} \label{pp}
x_{g_1}^{n_C}\cdot s=x_{g_2}^{n_C}\cdot s.\end{equation}
\end{lem}
\proof Equality (\ref{pp}) is proved in \cite{F-V}. The proof of
(\ref{oo}) is similar. \qed

\subsection{$C$-groups associated to equipped groups and the type homomorphism}
Let $(G,O)$ be an equipped group such that ${\bf{1}}\not\in O$ and
let the set $O$ be the union of $m$ conjugacy classes,
$O=C_1\cup\dots\cup C_m$.

A group $\hat G_O$, generated by an alphabet $Y_O=\{ y_g\mid g\in
O\}$ (so called $C$-generators) being subject to the relations
\begin{equation} \label{rel11} y_{g_1}
y_{g_2}=y_{g_2}y_{g_2^{-1}g_1g_2}=y_{g_1g_2g_1^{-1}}y_{g_1},\quad
y_{g_1},y_{g_2}\in Y_O,\end{equation} is called the {\it $C$-group
associated} to $(G,O)$. It is obvious that the maps $x_g\mapsto y_g$
and $y_g\mapsto g$ define two homomorphisms: $\beta : S(G,O)\to \hat
G_O$ and $\gamma :\hat G_O\to G$ such that
$\alpha=\gamma\circ\beta$. The elements of $\text{Im}\, \beta$ are
called the {\it positive} elements of $\hat G_O$.

A $C$-group $\hat G_O$, associated to an equipped group $(G,O)$, has
similar properties as the semigroup $S_O$ has. For example, like in
the case of factorization semigroups, it is easy to check that for
any $\hat g\in \hat G_O$ and any $g_1\in O$ the following relation
\begin{equation} \label{kkk} \hat g^{-1}y_{g_1}\hat g=y_{g^{-1}g_1g}\end{equation}
is a consequence of relations  (\ref{rel11}), where  $g=\gamma(\hat
g)$.

Denote by $\hat O$ the subset $\{ y_{g}\mid g\in O\}$ of $\hat G_O$.
It follows from relation (\ref{kkk}) that $\hat O$ is invariant
under inner automorphisms of $\hat G_O$.
\begin{claim} Let $(G,O)$ be an equipped group. Then the semigroups $S(G,O)$
and $S(\hat G_O,\hat O)$ are naturally isomorphic. \end{claim}
\proof By relations (\ref{kkk}), it is easy to see that the map
$\xi: S(\hat G_O,\hat O)\to S(G,O)$, given by $\xi(x_{y_g})=x_g$ for
$g\in O$,
is an isomorphism of semigroups. \qed \\

Applying relations (\ref{kkk}), it is easy to prove the following
proposition (see, for example, \cite{Kuz}).
\begin{prop} \label{simple1} For any equipped group $(G,O)$ we have $$Z(\hat G_O)=\gamma^{-1}(Z(G_O)),$$
where $Z(G_O)$ and $Z(\hat G_O)$ are the centers, respectively, of
$G_O$ and $\hat G_O$.\end{prop}

It is easy to see that the first homology group $H_1(\hat
G_O,\mathbb Z)= \hat G_O/[\hat G_O,\hat G_O]$ of $\hat G_O$ is a
free abelian group of rank $m$. Let $\text{ab}: \hat G_O\to H_1(\hat
G_O,\mathbb Z)$ be the natural epimorphism. The group $H_1(\hat
G_O,\mathbb Z)\simeq \mathbb Z^m$ is generated by
$\text{ab}(y_{g_i})=(0,\dots,0,1,0\dots,0)$ ($1$ stands on the
$i$-th place), where $g_i\in C_i$.

The homomorphism of semigroups $\tau = \text{ab}\circ \beta:
S(G,O)\to \mathbb Z_{\geq 0}^m\subset \mathbb Z^m$ is called the
{\it type homomorphism} and the image $\tau(s)$ of $s\in S(G,O)$ is
called the {\it type} of $s$. If $O$ consists of a single conjugacy
class, then the homomorphism $\tau$ can (and will) be identified
with the homomorphism $ln:S(G,O)\to \mathbb Z_{\geq 0}$.

\begin{lem} \label{simple2} Any element $\hat g$ of the $C$-group $\hat G_O$,
associated with an equipped group $(G,O)$, can be represented in the
form: \begin{equation} \label{diferen} \hat g=\hat g_1\hat
g_2^{-1},\end{equation} where $\hat g_1$, $\hat g_2$ are positive
elements. In particular, $\hat g\in \hat G_O'=[\hat G_O,\hat G_O]$
if and only if $ab(\hat g_1)=ab(\hat g_2)$ in representation {\rm
(\ref{diferen})} of $\hat g$ as a quotient of two positive elements
$\hat g_1$ and $\hat g_2$.
\end{lem}

\proof Write $\hat g$ in the form: $\hat
g=y_{g_{i_1}}^{\varepsilon_1}\dots y_{g_{i_k}}^{\varepsilon_k}$,
where $g_{i_j}\in O$ and $\varepsilon_j=\pm 1$. To prove lemma, it
suffices to note that by relations (\ref{rel11}) we have
$y_{g_2}^{-1}y_{g_1}=y_{g_2^{-1}g_1g_2}y_{g_2}^{-1}$ for any $g_1,
g_2\in O$.\qed

\begin{claim} Let $(G,O)$ be a finite equipped group. The homomorphism $\beta: S_O\to \hat G_O$ is an embedding
if and only if $O\subset Z(G_O)$, that is, $G_O$ is an abelian
group.
\end{claim}

\proof Let $O=C_1\cup \dots \cup C_m$ be the decomposition into the
union of conjugacy classes. %, ${\bf{1}}\not \in O$.
It is easy to
see that if $O\subset Z(G_O)$ then $\hat G_O\simeq \mathbb Z^{|O|}$,
where the isomorphism is induced by homomorphism $ab$, and in this
case the semigroup $S_O$ can be identified with semigroup $\mathbb
Z_{\geq 0}^{|O|}\subset \mathbb Z^{|O|}$.

If $O\not\subset Z(G_O)$, then there is $C_i\subset O$ consisting of
at least two elements, say $g_1$ and $g_2$. Let $n$ be their order
in $G$. Then it is easy to see that $x_{g_1}^n\neq x_{g_2}^n$ in
$S_O$. On the other hand, their images $y_{g_1}^n=\beta(x_{g_1}^n)$
and $y_{g_2}^n=\beta(x_{g_2}^n)$ coincide in $\hat G_O$. Indeed,
without loss of generality? we can assume that there is $g\in G_O$
such that $g_2=g^{-1}g_1g$. Consider an element $\hat g\in \gamma
^{-1}(g)$. Then
$$\hat g^{-1}y_{g_1}^n\hat g=(\hat g^{-1}y_{g_1}\hat
g)^n=y_{g^{-1}g_1g}^n=y_{g_2}^n,$$ but by Proposition \ref{simple1},
$y_{g_1}^n$ and $y_{g_2}^n$ belong to $Z(\hat G_O)$. Therefore
$y_{g_1}^n=y_{g_2}^n$. \qed

\subsection{Hurwitz equivalence.}\label{mu}  As above, let $O$ be a
subset of $G$ invariant under the action by inner automorphisms of
$G$. Consider the set
$$O^n=\{ (g_1,\dots ,g_n) \,\,\mid g_i\in O \}$$ of all ordered
$n$-tuples in $O$ and let $\text{Br}_n$ be the braid group with $n$
strings. We fix a set $\{ a_1,\dots ,a_{n-1} \} $ of so called {\it
standard {\rm (or} Artin{\rm )} generators} of $\text{Br}_n$, that
is, generators being subject to the relations \begin{equation}
\label{br}
\begin{array}{ll}
a_ia_{i+1}a_i & =a_{i+1}a_i a_{i+1} \qquad \qquad 1\leq i\leq n-1 ,  \\
a_ia_{k} & =a_{k}a_i  \qquad \qquad \qquad \, \, \mid i-k\mid \,
\geq 2.
\end{array}
\end{equation}
The group $\text{Br}_n$ acts on $O^n$ as follows
$$((g_1,\dots ,g_{i-1},g_i,g_{i+1},g_{i+2}, \dots ,g_n))a_i=
(g_1,\dots ,g_{i-1},g_ig_{i+1}g_i^{-1},g_{i},g_{i+2}, \dots
,g_n)).$$ Usually, the action of the standard generators $a_i\in
\text{Br}_n$ and their inverses on $O^n$ is  called {\it Hurwitz
moves}. Two elements in $O^n$ are called  {\it Hurwitz equivalent}
if one can be obtained from the other by a finite sequence of
Hurwitz moves, that is, if they belong to the same orbit under the
action of $\text{Br}_n$.

There is a natural map ({\it product map}) $\alpha: O^n\to G$
defined by
$$\alpha((g_1,\dots ,g_n))=g_1\dots g_n$$
and an element $(g_1,\dots ,g_n)\in O^n$ is called a { \it
factorization of $g=\alpha((g_1,\dots ,g_n))\in G$ with factors in}
$O$.

There is a natural map $\varphi:O^n\to S(G,O)$ sending $(g_1,\dots
,g_n)$ to $s=x_{g_1}\cdot\, .\, .\, .\, \cdot x_{g_n}$.

\begin{claim} \label{cl1.3} Two factorizations $y$ and
$z\in O^n$ are Hurwitz equivalent if and only if $\varphi
(y)=\varphi (z)$.
\end{claim}

\proof Evident. \qed

\begin{rem} {\rm In what follows, according with Claim \ref{cl1.3}, we
identify classes of Hurwitz equivalent factorizations in $O$ with
their images in $S(G,O)$.}
\end{rem}

Define also the {\it conjugation action} of $G$ on $O^n$:
$$ \lambda(g)((g_1,\dots ,g_n))= (g^{-1}g_1g,\dots ,g^{-1}g_ng).$$
Obviously, this action is compatible under the map $\varphi$ with
the conjugation action of $G$ on $S(G,O)$ defined above.

Denote by $W=W(O)$ the set of words in the alphabet $X=X_{O\setminus
\{\bf{1}\}}$ and by $W_n$ its subset consisting of the words of
length $n$. In what follows, the elements of the set $O^n$ will be
identified with the elements of $W_n$ (identification: $(g_1,\dots
,g_n)\in O^n\leftrightarrow x_{g_1}\dots x_{g_n}\in W_n$) and we put
$$W(s)=\{ w\in W\mid \varphi (w)=s\in S(G,O)\}.$$

\subsection{Finite presentability of some subsemigroups of $S(G,O)$}
Let $(G,O)$ be a finite equipped group. By definition, the semigroup
$S_O$ is finitely presented. From geometric point of view the most
interesting subsemigroups of $S_G$ are $S_{O,\bf{1}}$ and
$S_{O,{\bf{1}}}^G=\{ s\in S_{O,{\bf{1}}}\mid G_s=G\}$. (Note that
$S_{O,{\bf{1}}}^G$ is non-empty if and only if $G_O=G$.) In this
subsection, we will show that the semigroups $S_{O,\bf{1}}$ are
finitely presented, but for the semigroups $S_{O,{\bf{1}}}^G$ the
property to be finitely presented (and, moreover, to be finitely
generated) is not obligatory.

Let $N=|G|$ be the order of $G$ and $\mathcal C=\{ C_1,\dots, C_m\}$
be the set of conjugacy classes of $G$ such that $O=\cup C_i$. For
$C\in\mathcal C$ let $n_{C}=n_g$ be the order of $g\in C$. In each
$C\in\mathcal C$ we choose and fix an element $g_C\in C$.

It is evident that a necessary condition for a subsemigroup $S$ of
$S_O$ to be finitely  generated is that its image $\tau(S)$ is a
finitely generated  semigroup, where $\tau :S_O\to \mathbb Z_{\geq
0}^m$ is the type homomorphism.
\begin{thm} \label{fpr} A factorization semigroup $S_{O,\bf{1}}$ over a finite
group $G$ is finitely presented.
\end{thm}
\proof Let $O=C_1\cup \dots \cup C_m$ be the decomposition into the
union of conjugacy classes and let ${\bf{1}}\not\in O$. We numerate
the elements of $O=\{ g_1,\dots, g_{K}\}$ so that $g_i=g_{C_i}$ for
$i=1,\dots,m$.

For any $g\in O$ we have $s_{g}=x_g^{n_g}\in S_{O,\bf{1}}$. Let
$F=\{ s_1, \dots ,s_M\}$ be the set of elements of $S_{O,\bf{1}}$ of
length less or equal to $K^N$, where $N=|G|$, and we assume also
that $s_i=s_{g_i}=x_{g_{i}}^{n_{g_i}}$ for $i\leq K$. Let us show
that the elements $s_1,\dots, s_M\in F$ generate the semigroup
$S_{O,\bf{1}}$.
\begin{lem}\label{si} An element  $s\in S_{O,\bf{1}}$ of
length $ln(s)> K^N$ can be written in the following form:
$$s=s_{i_1}^{n_1}\cdot \, .\, .\, .\, \cdot s_{i_l}^{n_l}\cdot
\overline s,$$ where $1\leq i_1\leq \dots\leq i_l\leq K$  and
$\overline s\in S_{O,\bf{1}}$ with $ln(\overline s)\leq K^N$.
\end{lem}
\proof If $ln(s)>K^N$, then in a presentation of $s$ as a product
$x_{g_1}\cdot\, .\, .\, .\, \cdot x_{g_{ln(s)}}$ there are at least
$N$ coinciding factors $x_g$ for some $g\in O$. Since $n_g\leq N$,
moving $n_g$ of these factors to the left (by means of relations
(\ref{rel1})), we obtain that $s=s_g\cdot s'$, where $s'\in
S_{O,\bf{1}}$ and $ln(s')<ln(s)$. \qed
\\

It follows from Lemma \ref{si} that $S_{O,\bf{1}}$ is generated by
the elements $s\in S_{O,\bf{1}}$ of length $ln(s)\leq K^N$, that is,
$S_{O,\bf{1}}$ is finitely generated.

To show that $S_{O,\bf{1}}$ is finitely presented, let us divide the
set of all relations as follows. The first set $R_1$ of relations
consists of relations:
$$ s_i\cdot s_j=s_j\cdot s_i, \quad \, %\text{for}\,
s_i, s_j\in F. %1\leq i<j\leq M.
$$

Denote by ${\bf{k}}=(k_1,\dots,k_M)$ an ordered collection of
non-negative integers and put $s_{\bf{k}}= s_1^{k_{1}}\cdot\, .\,
.\, .\cdot s_M^{k_{M}}$. In view of the existence of relations
$R_1$, we can assume that all other relations in $S_{O,\bf{1}}$
connecting the generators $s_1,\dots, s_M$ have the following form:
\begin{equation}
\label{rela} s_{{\bf{k}}_1}= s_{{\bf{k}}_2}.\end{equation} Note that
if we have a relation of form (\ref{rela}), then
$G_{s_{{\bf{k}}_1}}= G_{s_{{\bf{k}}_2}}$ and $\tau (s_{{\bf{k}}_1})=
\tau (s_{{\bf{k}}_2})$.

Consider the set $\overline R_2$ of all relations of form
(\ref{rela}) for which $G_{s_{{\bf{k}}_1}}$ is a proper subgroup of
$G$. By induction, we can assume that the semigroups
$S(\Gamma,\overline O)_{\bf{1}}$ are finitely presented for all
equipped groups $(\Gamma,\overline O)$ of order less than $N$. Since
there are only finitely many proper subgroups of $G$ and the
embeddings $(G_{s_{{\bf{k}}_1}},O\cap
G_{s_{{\bf{k}}_1}})\hookrightarrow (G,O)$ define the embeddings
$S(G_{s_{{\bf{k}}_1}},O\cap
G_{s_{{\bf{k}}_1}})_{\bf{1}}\hookrightarrow S_{O,\bf{1}}$, we obtain
that there is a finite set of relations $R_2\subset \overline R_2$
generating all relations of $\overline R_2$.

Denote by $R_3$ the set of all relations in $S_{O,\bf{1}}$ of the
form $s_{\bf{k}_1}=s_{\bf{k}_2}$ which are not contained in $R_1\cap
R_2$ and such that $ln(s_{\bf{k}_1})\leq K^N$.  It is easy to see
that $R_3$ is a finite set.

For each element $s_i$ of the set of generators of $S_{O,\bf{1}}$
with $i\geq K+1$, we put $$n_i=\min_n\{ ln(s_i^n)>K^N\}-1.$$ From
Lemma \ref{si} it follows
\begin{lem} \label{yyy} For any $i\geq K+1$
the element $s_i^{n_i+1}$ can be written in the following form:
\begin{equation} \label{mmm} s_i^{n_i+1}=(\prod_{j=1}^{K}s_j^{a_j})\cdot s_l,\end{equation} where
${\bf{a}}=(a_1,\dots,a_{K}) $ is a collection of non-negative
integers and $s_l\in F$ is a generator with index  $l\geq K+1$.
\end{lem}

Denote by $R_4$ the set of relations of form (\ref{mmm}). It is a
finite set. By Lemma \ref{yyy}, applying relations of the set
$R_1\cup R_4$, each element $s\in S_{O,\bf{1}}$ can be written in
the form: $ s=s_{\bf{k}}$, where ${\bf{k}}=(k_1,\dots, k_M)$
satisfies the following condition: $k_i\leq n_i$ for $i\geq K+1$.

An element $s_{\bf{k}}$ is called {\it $\Gamma$-primitive} if in
${\bf{k}}=(k_1,\dots, k_M)$ all $k_i\leq 1$ for $i\leq K$, $k_i\leq
n_i$ for $i\geq K+1$, and $G_{s_{\bf{k}}}=\Gamma$. By Lemma
\ref{fried}, for each $G$-primitive element $s_{\bf{k}}$ we have the
following relations in $S_{O,\bf{1}}$:
$$s_{i}\cdot s_{\bf{k}}=s_j\cdot s_{\bf{k}},$$ where
$i\leq m$ and $j\leq K$ is such that $g_j\in C_i$. Denote by $R_5$
the set of all such relations. Obviously, $R_5$ is a finite set.

Let $s\in S_{O,\bf{1}}$ be such that $G_s=G$.   Applying relations
of $R_5$, as above it is easy to show that $s$ can be written in the
form: \begin{equation} \label{new} s=(\prod_{j=1}^{m}s_j^{a_j})\cdot
s_{\bf{k}},\end{equation} where $s_{\bf{k}}$ is some $G$-primitive
element. Denote by $\overline R_6$ the set of relations in
$S_{O,\bf{1}}$ of the form:
\begin{equation} \label{vvv} (\prod_{j=1}^{m}s_j^{b_{j,1}})\cdot
s_{{\bf{k}}_1}=(\prod_{j=1}^{m}s_j^{b_{j,2}})\cdot s_{{\bf{k}}_2},
\end{equation} where $s_{{\bf{k}}_1}$ and $s_{{\bf{k}}_2}$ are $G$-primitive elements.

To complete the proof of Theorem \ref{fpr}, it suffices to show that
the relations of $\overline R_6$ are consequences of a finite set of
relations $R_6$. Since there are only finitely many $G$-primitive
elements, it is suffices to show that for fixed $G$-primitive
elements $s_{{\bf{k}}_1}$ and $s_{{\bf{k}}_2}$ relations of form
(\ref{vvv}) are consequences of a finite set of relations. For this
purpose, consider the semigroup $\mathbb Z_{\geq 0}^m=\{ {\bf{a}}=
(a_1,\dots,a_m)\in \mathbb Z^m\mid a_i\geq 0\}$.

A subsemigroup $S$ of $\mathbb Z_{\geq 0}^m$ is called {\it
non-perforated} if for any ${\bf{a}}\in S$ and any ${\bf{b}}\in
\mathbb Z_{\geq 0}^m$ the element ${\bf{a}+\bf{b}}\in S$. Note that
if $S_1$ and $S_2$ are non-perforated subsemigroups, then $S_1\cup
S_2$ and $S_1\cap S_2$ are also non-perforated subsemigroups. An
element ${\bf{a}}$ of a non-perforated subsemigroup $S$ is called an
{\it origin} of $S$ if there does not exist elements ${\bf{b}}\in S$
and ${\bf{c}}\in \mathbb Z_{\geq 0}^m\setminus \{ {\bf{0}}\}$ such
that ${\bf{a}=\bf{b}+\bf{c}}$. Denote by $O(S)$ the set of origins
of a non-perforated subsemigroup $S$. A non-perforated subsemigroup
$S$ with a single origin is called {\it prime}. It is easy to see
that if ${\bf{a}}$ is the origin of a prime non-perforated
subsemigroup $S$, then
$$S=F_{\bf{a}}=\{ {\bf{c}}={\bf{a}}+{\bf{b}}\in \mathbb Z_{\geq
0}^m\mid {\bf{b}}\in \mathbb Z_{\geq 0}^m\} .$$ It is obvious that a
non-perforated subsemigroup  $S$ can be represented as a union of
prime non-perforated subsemigroups, for example,
$$S=\bigcup_{{\bf{a}}\in S}F_{\bf{a}}.$$
Let $A$ be a subset of $S$ and let $S$ be represented as the union
of prime non-perforated subsemigroups,
\begin{equation} \label{rrrrr} S=\bigcup_{{\bf{a}}\in
A}F_{\bf{a}}. \end{equation} We say that representation
(\ref{rrrrr}) is {\it minimal} if
$$S\neq \bigcup_{{\bf{a}}\in A\setminus \{
{\bf{a}}_0\} }F_{\bf{a}} $$ for any ${\bf{a}}_0\in A$.

\begin{claim} \label{hhh} For a non-perforated subsemigroup
$S\subset \mathbb Z_{\geq 0}^m$ there is the unique minimal
representation as the union of prime non-perforated subsemigroups,
namely, $$S=\bigcup_{{\bf{a}}\in O(S)}F_{\bf{a}}.$$
\end{claim}
\proof It follows from the definition of origins that if $S=\cup
F_{{\bf{a}}_i}$ is a representation as the union of prime
non-performed subsemigroups and $\bf{a}$ is an origin of $S$, then
${\bf{a}}={\bf{a}_i}$ for some $i$.

Assume that $$C=S\setminus\bigcup_{{\bf{a}}\in O(S)}F_{\bf{a}}$$ is
not empty, then there is ${\bf{c}}_0=(c_{1,0},\dots, c_{m,0})\in C$
such that $c_{m,0}=\min c_m$ for $(c_1,\dots, c_m)\in C$,
$c_{m-1,0}=\min c_{m-1}$ for $(c_1,\dots,c_{m-1}, c_{m,0})\in C$,
$\dots$, $c_{1,0}=\min c_{1}$ for $(c_1,c_{2,0}\dots, c_{m,0})\in
C$. It is obvious that ${\bf{c}}_0$ is an origin of $S$. \qed

\begin{prop}\label{n-p} Every increasing sequence of non-perforated
subsemigroups of $\mathbb Z_{\geq 0}^m$,
$$S_1\subset S_2\subset S_3\subset \dots ,$$
such that $S_i\neq S_{i+1}$ is finite. \end{prop}

\proof Proposition is obvious if $m=1$. let us use the induction on
$m$. Consider an increasing sequence of non-perforated subsemigroups
$S_1\subset S_2\subset S_3\subset \dots \subset\mathbb Z_{\geq
0}^m$, $m\geq 2$. Denote by $P_j=\{ (z_1,\dots,z_m)\in \mathbb
Z_{\geq 0}^m\mid z_m=j\}$, and $S_{i,j}=S_i\cap P_j$. Then $S_{i,j}$
also can be considered as a non-perforated subsemigroup of $\mathbb
Z_{\geq 0}^{m-1}$ (if we forget about the last coordinate). By
inductive assumption, increasing sequences $S_{1,j}\subset
S_{2,j}\subset S_{3,j}\subset \dots$ must stop for each $j$. Denote
by $\overline S_j=S_{i(j),j}$ the first biggest semigroups in these
sequences.

Consider a map $sh: \mathbb Z_{\geq 0}^m\to \mathbb Z_{\geq 0}^m$ is
given by
$$sh((z_1,\dots,z_{m-1},z_m))=(z_1,\dots,z_{m-1},z_m+1).$$  It follows from definition of non-perforated
subsemigroups that $sh:S_{i,j}\to S_{i,j+1}$ is an embedding of
semigroups. Therefore we can (and will) identify a semigroup
$S_{i,j}$ with subsemigroup $sh^n(S_{i,j})$ of $S_{i,j+n}$. It
follows from definition of non-performed subsemigroups that if $j_1
< j_2$, then $\overline S_{j_1}=S_{i(j_1),j_1}\subset \overline
S_{j_2}=S_{i(j_2),j_2}$. As a result we obtain an increasing
sequence of non-perforated subsemigroups  $$S_{i(0),0}\subset
S_{i(1),1}\subset S_{i(2),2}\subset \dots \subset\mathbb Z_{\geq
0}^{m-1}.$$ It must stop. It is easy to see that if $S_{i(j_0),j_0}$
is the biggest semigroup, then
$S_{i(j_0)}=S_{i(j_0)+1}=S_{i(j_0)+2}=\dots $. \qed

\begin{cor} \label{non-per} The set of origins $O(S)$ of a non-perforated subsemigroup
$S\subset \mathbb Z_{\geq 0}^m$ is  non-empty and finite.
\end{cor}

\proof If the set $O(S)=\{ {\bf{a}}_1,{\bf{a}}_2,{\bf{a}}_3,\dots
\}$ is infinite, then by Claim \ref{hhh} we will have an infinite
increasing sequence $$F_{{\bf{a}}_1}\subset F_{{\bf{a}}_1}\cup
F_{{\bf{a}}_2}\subset F_{{\bf{a}}_1}\cup F_{{\bf{a}}_2}\cup
F_{{\bf{a}}_3}\subset \dots \,\, .$$ which contradicts Proposition
\ref{n-p}. \qed \\

Let us return to the proof that the relations of the set $\overline
R_6$ are consequences of a finite set of relations $R_6$. For this
purpose, note that if
\begin{equation} \label{rr}
(\prod_{j=1}^{m}s_j^{b_{j,1}})\cdot
s_{{\bf{k}}_1}=(\prod_{j=1}^{m}s_j^{b_{j,2}})\cdot
s_{{\bf{k}}_2}\end{equation} is a relation, then
$$(b_{1,1}n_{C_1},\dots,b_{m,1}n_{C_m})+\tau(s_{{\bf{k}}_1})=
(b_{1,2}n_{C_1},\dots,b_{m,2}n_{C_m})+\tau(s_{{\bf{k}}_2}).$$
Therefore if
$\tau(s_{{\bf{k}}_j})=(\alpha_{1,j},\dots,\alpha_{m,j})$, then
$\alpha_{i,1}\equiv \alpha_{i,2}( \text{mod}\, n_{C_i})$ for all
$i$. Put $a_{i,1,0}=b_{i,1}-b_{i,2}$ if $\alpha_{i,2}\geq
\alpha_{i,1}$ and $a_{i,1,0}=0$ if otherwise. Respectively, put
$a_{i,2,0}=b_{i,2}-b_{i,1}$ if $\alpha_{i,1}\geq \alpha_{i,2}$ and
$a_{i,2,0}=0$ if otherwise. We have
$$n_{C_i}a_{i,1,0}+\alpha_{i,1}=n_{C_i}a_{i,2,0}+\alpha_{i,2}$$
and $a_{i,1,0}$, $a_{i,2,0}$ are defined uniquely by $\alpha_{i,1}$,
$\alpha_{i,2}$, and $n_{C_i}$. Moreover, if we denote
$a_{i,j}=b_{i,j}-a_{i,j,0}$, then $a_{i,1}=a_{i,2}\geq 0$ for
$i=1,\dots,m$, and each relation of the form (\ref{rr}) can be
rewritten in the form
\begin{equation} \label{rrr} (\prod_{j=1}^{m}s_j^{a_{j}})\cdot
(\prod_{j=1}^{m}s_j^{a_{j,1,0}})\cdot
s_{{\bf{k}}_1}=(\prod_{j=1}^{m}s_j^{a_{j}})\cdot
(\prod_{j=1}^{m}s_j^{a_{j,2,0}})\cdot s_{{\bf{k}}_2}, \end{equation}
where $a_j=a_{j,1}=a_{j,2}$.

If (\ref{rrr}) is a relation in $S_{O,\bf{1}}$, then
$$(\prod_{j=1}^{m}s_j^{a_{j}+b_j})\cdot
(\prod_{j=1}^{m}s_j^{a_{j,1,0}})\cdot
s_{{\bf{k}}_1}=(\prod_{j=1}^{m}s_j^{a_{j}+b_j})\cdot
(\prod_{j=1}^{m}s_j^{a_{j,2,0}})\cdot s_{{\bf{k}}_2}$$ is also a
relation for each ${\bf{b}}=(b_1,\dots, b_m)\in \mathbb Z_{\geq
0}^m$ and it is a consequence of relation (\ref{rrr}).

It follows from consideration above that the set $\{
(a_1,\dots,a_m)\}$ of exponents interning into the relations written
in the form (\ref{rrr}) for fixed $s_{{\bf{k}}_1}$ and
$s_{{\bf{k}}_2}$ form a non-perforated subsemigroup
$F_{s_{{\bf{k}}_1},s_{{\bf{k}}_2}}$ of $\mathbb Z_{\geq 0}^m$. By
Corollary \ref{non-per}, the set
$O(F_{s_{{\bf{k}}_1},s_{{\bf{k}}_2}})$ of its origins is finite. It
is easy to see that the relations (\ref{rrr}) for fixed
$s_{{\bf{k}}_1}$ and $s_{{\bf{k}}_2}$ are consequences of the
relations corresponding to the origins of
$F_{s_{{\bf{k}}_1},s_{{\bf{k}}_2}}$, and since there are only
finitely many $G$-primitive elements, we obtain that the relations
of $\overline R_6$ are consequences of some finite subset $R_6$ of
$\overline R_6$.

To complete the proof of Theorem \ref{fpr}, it suffices to note that
all relations are consequences of the relations belonging to
$R_1\cup \dots \cup R_6$ which is a finite set. \qed \\

Note that not any subsemigroup  $S_{O,{\bf{1}}}^G$ of $S_G$ is
finitely generated. For example, let $G\simeq (\mathbb Z/2\mathbb
Z)^2$ be generated by two elements $g_1$ and $g_2$. If $O=\{ g_1,
g_2\}$, then $S_{O,{\bf{1}}}^G$ is isomorphic to the semigroup
$$S=\{ (a_1,a_2)\in \mathbb Z_{\geq 0}^2\mid a_1>0, \, a_2>0\}$$
which is not finitely generated.

\begin{prop} \label{non-gen} Let $(G,O)$ be a finite equipped group.
Assume that $O$ is the union of conjugacy classes, $O=C_1,\cup \dots
\cup C_m$, such that for each $i$ the elements of $C_i$ generate the
group $G$. Then the subsemigroup $S_{O,{\bf{1}}}^G$ of $S_G$ is
finitely presented.
\end{prop}

\proof In notations used in the proof of Theorem \ref{fpr}, denote
$$s_{C_i}=\prod_{g_l\in C_i} x_{g_l}^{n_{C_i}}=\prod_{g_l\in C_i}
s_{l}.$$ We have $s_{C_i}\in S_{O,{\bf{1}}}^G$, since the elements
$g_l\in C_i$ generate $G$.

As it was shown in the proof of Theorem \ref{fpr}, any element $s\in
S_{O,{\bf{1}}}^G$ can be written in the form (\ref{new}):
$$s=(\prod_{i=1}^ms_{i}^{a_i})\cdot s_{\bf{k}},$$ where
$s_{\bf{k}}$ is some $G$-primitive element of $S_{O,{\bf{1}}}^G$. If
$a_i\geq |C_i|$, then by Lemma \ref{fried},
$$s_{i}^{a_i}\cdot s_{\bf{k}}= s_{C_i}\cdot s_{i}^{a_i-|C_i|}\cdot s_{\bf{k}}.$$
Therefore any element $s\in S_{O,{\bf{1}}}^G$ can be written in the
form
\begin{equation} \label{ww} s= (\prod_{i=1}^m
s_{C_i}^{b_i})\cdot (\prod_{i=1}^ms_{i}^{a_i})\cdot
s_{\bf{k}},\end{equation} where $(b_1,\dots, b_k)\in \mathbb Z_{\geq
0}^k$ and $0\leq a_i<|C_i|$, and $s_{\bf{k}}$ is a $G$-primitive
element. Since there are only finitely many expressions of the form
\begin{equation} \label{www} (\prod_{i=1}^ms_{i}^{a_i})\cdot s_{\bf{k}},\end{equation}
where $0\leq a_i<|C_i|$, and $s_{\bf{k}}$ is a $G$-primitive
element, the end of the proof of Proposition \ref{non-gen} coincides
with the proof of Theorem \ref{fpr}. \qed

\subsection{Stabilizing elements} If $G$ is an abelian finite group, then it is obvious that
the type homomorphism $\tau :S_G\to \mathbb Z_{\geq 0}^{|G|-1}$ is
an isomorphism. If $G$ is not an abelian group and $c(G)$ is the
number of conjugacy classes of its elements $g\neq \bf{1}$, then the
type homomorphism $\tau :S_G\to \mathbb Z_{\geq 0}^{c(G)}$ is a
surjective, but not injective homomorphism, and one of the main
problems is to describe the preimages $\tau^{-1}({\bf{a}})$ of
elements ${\bf{a}}\in \mathbb Z_{\geq 0}^{c(G)}$ (in particular, to
describe the set of elements ${\bf{a}}\in \mathbb Z_{\geq 0}^{c(G)}$
for which each element $s\in \tau^{-1}({\bf{a}})$ is uniquely
determined by their value $\alpha(s)\in G$).

\begin{prop} \label{f} Let $S_{O,{\bf{1}}}^G$ be as in Proposition
{\rm \ref{non-gen}}. Then there is a constant $c=c(G,O)$ such that
for any ${\bf{a}}\in \mathbb Z_{\geq 0}^m$ the number
$|\tau^{-1}({\bf{a}})|$ of preimages of $\bf{a}$ under the
homomorphism $\tau: S_{O,{\bf{1}}}^G\to \mathbb Z_{\geq 0}^m$ is
less than $c$. \end{prop}

\proof In the proof of Proposition \ref{non-gen} it was shown that
any element $s\in S_{O,{\bf{1}}}^G$ can be written in the form
(\ref{ww}). Therefore Proposition \ref{f} follows from that the
number of different expressions of the form (\ref{www}) is finite.
\qed \\

Note that Proposition \ref{f} is false if we consider the semigroup
$S_{O,{\bf{1}}}$  instead of $S_{O,{\bf{1}}}^G$, see, for example,
Corollary \ref{uniq}.

An element $s\in S(G,O)$ is called  {\it stabilizing} if $s\cdot
s_1=s\cdot s_2$ for any $s_1,s_2\in S(G,O)$ such that
$\tau(s_1)=\tau(s_2)$ and $\alpha(s_1)=\alpha(s_2)$. A semigroup
$S(G,O)$ is called {\it stable} if it possesses a stabilizing
element.
\begin{claim} If $s$ is a stabilizing element of $S(G,O)$, then for any $s_1\in S(G,O)$ the element
$s\cdot s_1$ is also a stabilizing element. In particular, if
$S(G,O)$ is a stable semigroup, then there is a stabilizing element
$s\in S(G,O)$ such that $\alpha(s)={\bf{1}}$.
\end{claim}
\proof Evident. \qed \\

Conway -- Parker Theorem (see  Appendix in \cite{F-V}) gives some
sufficient condition for a semigroup $S_G$ to be stable. To
formulate this theorem, recall that a {\it Schur covering group} $R$
of a finite group $G$ is a group of maximal order with the property
that $R$ has a subgroup $M\subset R'\cap Z(R)$ satisfying $R/M\simeq
G$, where $R'=[R,R]$ is the commutator subgroup and $Z(R)$ is the
center of $R$. Such an $R$ always exists (but non necessarily
unique). The group $M$ isomorphic to the Schur multiplier
$M(G)=H^2(G,\mathbb C^*)$ of $G$. The Schur multiplier $M(G)$ is
said to be {\it generated by commutators} if $M\cap
\{g^{-1}h^{-1}gh\mid g,h\in R\}$ generates $M$.
\begin{thm} {\rm (Conway -- Parker) (\cite{F-V})} \label{C-P} Let $G$ be a finite group,
$O=G\setminus {\bf{1}}=C_i\cup\dots\cup C_m$ the decomposition into
the union of conjugacy classes, and denote
$$\overline s=\prod_{g\in G\setminus \{ {\bf{1}}\}}x_g^{n_g}\in S_{G},$$
where $n_g$ is the order of $g$ in $G$. Assume that the Schur
multiplier $M(G)$ of $G$ is generated by commutators. Then there is
a constant $n=n(G)$ such that $\overline s^n$ is a stabilizing
element of $S_G$.
\end{thm}

Note that a Schur covering group $G$ of a finite group $H$ satisfies
the conditions of Conway -- Parker Theorem (see \cite{F-V}).

In the next section we will prove that factorization semigroups
$S_{\mathcal S_d}$ over symmetric groups $\mathcal S_d$ are also
stable. On the other hand, there are many finite equipped groups
$(G,O)$ for which
$S(G,O)$ is not a stable semigroup. %such that
%for almost all $t %${\bf{a}}
%$\in \tau(S(G,O))$ the elements of $\tau^{-1}(t)$ %({\bf{a}})$
%are not uniquely defined by the value $\alpha(s)\in G$.

\begin{prop} \label{non-one} Let $(H, \tilde O)$ be a finite equipped group such that
\begin{itemize}
\item[$(i)$] the elements of $\tilde O$ generate the group $H$;
\item[$(ii)$] $H'\cap Z(H)\neq {\bf{1}}$;
\item[$(iii)$] $\tilde g_1\tilde g_2^{-1}\not\in Z(H)\setminus \{ {\bf{1}} \}$
for all $\tilde g_1, \tilde g_2\in \tilde O$.
\end{itemize}
Let $f: H\to H/Z(H)=G$ be the natural epimorphism and put
$O=f(\tilde O)\subset G$. Then there are at least two elements
$s_{1}, s_{2}\in S_{O,{\bf{1}}}^G$ such that $\tau(s\cdot
s_{1})=\tau(s\cdot s_{2})$, but $s\cdot s_{1}\neq s\cdot s_{2}$ for
all  $s \in S_{O,{\bf{1}}}^G$.

In particular, if $\tilde O$ consists of a single conjugacy class of
$H$, then there is a constant $N\in \mathbb N$ such that for any
$t\in \tau(S_{O,{\bf{1}}}^G)\cap \mathbb Z_{\geq N}$ there are at
least two elements $s_1,s_2\in S_{O,{\bf{1}}}^G$  such that
$\tau(s_1)=\tau(s_2)=t$, but $s_1\neq s_2$.
\end{prop}

\proof By $(i)$, the elements of $O$ generate the group $G$. By
$(iii)$, the surjective map $f_{|\tilde O} : \tilde O\to O$ is a
bijection, and if we denote $g_i=f(\tilde g_i)$ for $\tilde g_i\in
\tilde O$, then the equality $g_i^{-1}g_jg_i=g_k$ holds in $G$ for
some $g_i,g_j,g_k\in O$ if and only if the equality $\tilde
g_i^{-1}\tilde g_j\tilde g_i=\tilde g_k$ holds in $H$. Therefore the
induced homomorphism $f_* : S_{\tilde O}\to S_O$ (sending the
generators $x_{\tilde g_i}$ of $S_{\tilde O}$ to the generators
$x_{g_i}$ of $S_O$) is an isomorphism of semigroups. In particular,
the restriction of $f_*$ to $S_{\tilde O,Z(H)}^H=\{ \tilde s\in
S_{O}^H\mid \alpha(\tilde s)\in Z(H)\}$ gives an isomorphism between
$S_{\tilde O,Z(H)}^H$ and $S_{O,{\bf{1}}}^G$. In addition, the
homomorphism $f$ induces a surjective homomorphism $f_*:\hat
H_{\tilde O}\to \hat G_O$ of $C$-groups associated to $(H,\tilde O)$
and $(G,O)$ (sending the generators $y_{\tilde g_i}$ of $\hat
H_{\tilde O}$ to the generators $y_{g_i}$ of $\hat G_O$) such that
the following diagram \\

\begin{picture}(-30,0)(-30,0)
\put(120,-5){$ S_{\tilde O}  \buildrel{\beta}\over\longrightarrow
\hat H_{\tilde O} \buildrel{\gamma}\over\longrightarrow  H $}

\put(120,-50){$S_O \buildrel{\beta}\over\longrightarrow \hat G_O
\buildrel{\gamma}\over\longrightarrow G $}
\put(126,-10){\vector(0,-1){25}} \put(116,-23){$f_*$}
\put(128,-23){$\simeq$} \put(166,-10){\vector(0,-1){25}}
\put(168,-23){$f_*$} \put(208,-10){\vector(0,-1){25}}
\put(210,-23){$f$} %\put(310,-23){$(*)$}

\end{picture} \vspace{2cm} \newline
is commutative and  such that the induced homomorphism
$$f_{**}:H_1(\hat H_{\tilde O},\mathbb Z)\to H_1(\hat G_{O},\mathbb
Z)$$ is an isomorphism compatible with isomorphism $f_* : S_{\tilde
O}\to S_O$ (that is, if $s=f_*(\tilde s)$, then
$\tau(s)=f_{**}(\tau(\tilde s))$.  Therefore to prove the first part
of Proposition \ref{non-one}, it suffices to show that there are two
elements $\tilde s_1, \tilde s_2\in S_{\tilde O,Z(H)}^H$ such that
$\tau(\tilde s_1)=\tau(\tilde s_2)$, but $\alpha(\tilde s_1)\neq
\alpha(\tilde s_2)$. Indeed, for such two elements we will have that
$\tau(\tilde s\cdot \tilde s_1)=\tau(\tilde s\cdot \tilde s_2)$, but
$\alpha(\tilde s\cdot \tilde s_1)\neq \alpha(\tilde s\cdot \tilde
s_2)$ for all $\tilde s\in  S_{\tilde O,Z(H)}^H$. Therefore
$s_1=f_*(\tilde s_1)$ and $s_2=f_*(\tilde s_2)$ are non-equal
elements of $S_{O,{\bf{1}}}$ and $\tau(s\cdot s_1)=\tau(s\cdot
s_2)$, but $s\cdot s_1\neq s\cdot s_2$ for all elements $s\in
S_{O,{\bf{1}}}^G$ in view of isomorphism $f_*:S_{\tilde
O,Z(H)}^H\buildrel{\simeq}\over\longrightarrow S_{O,{\bf{1}}}^G$.

It follows from Proposition \ref{simple1} that for any subgroup
$\hat H_1$ of $\hat H_{\tilde O}$ we have $$\gamma(\hat H_1\cap
Z(\hat H_{\tilde O})) =\gamma(\hat H_1)\cap Z(H),$$ in particular,
$$ \gamma(\hat H_{\tilde O}'\cap Z(\hat H_{\tilde O}))=H'\cap Z(H).$$
Hence, by condition $(ii)$, there is an element $\hat h\in \hat
H_{\tilde O}'\cap Z(\hat H_{\tilde O})\setminus \{ {\bf{1}}\}$. By
Lemma \ref{simple2}, $\hat h=\hat h_1\hat h_2^{-1}$, where $\hat
h_1=\beta (\hat s_1)$ and $\hat h_2=\beta(\hat s_2)$ for some $\hat
s_1,\hat s_2\in S_{\hat O}$ (that is, $\hat h_1$ and $\hat h_2$ are
positive elements). Since $\hat h\in \hat H_{\tilde O}'$, we have
$ab(\hat h_1)=ab(\hat h_2)$.

Each element of a finite group $H$ can be expressed as a positive
word in its generators. Therefore, by condition $(i)$, there are
$\hat s\in S_{\tilde O}$ and the positive element $\hat g=\beta(\hat
s)\in \hat H_{\tilde O}$ such that $\gamma(\hat g)=\gamma(\hat
h_2^{-1})$. Denote also by $\hat s_0=\prod_{\tilde g_i\in \tilde
O}x_{\tilde g_i}^{n_i}\in S_{\tilde O,{\bf{1}}}^H$, where $n_i$ is
the order of $\tilde g_i$. Then $\tilde s_1=\hat s_0\cdot \hat
s\cdot \hat s_1$ and $\tilde s_2=\hat s_0\cdot \hat s\cdot \hat s_2$
are desired elements.

To prove the second part of Proposition \ref{non-one}, let us choose
elements $\overline s_1,\dots ,\overline s_n$ generating
$S_{O,{\bf{1}}}^G$ (by Proposition \ref{non-gen}, the semigroup
$S_{O,{\bf{1}}}^G$ is finitely generated in the case when $O$
consists of a single conjugacy class) and let $s_1$, $s_2$ be
elements the existence of which was proved in the first part of the
proof. Denote by $t_0=\tau(s_1)=\tau(s_2)$ and $t_i=\tau(\overline
s_i)$, $i=1,\dots,n$, and let $GCD(t_1,\dots, t_n)=d$, $t_i=a_id$.
Then the type $\tau(s)$ of any element of $S_{O,{\bf{1}}}^G$ is
divisible by $d$. Let us show that there is a constant $M\in \mathbb
N$ such that for any $j\in \mathbb N$ there is an element $s\in
S_{O,{\bf{1}}}^G$ with $\tau(s)=(M+j)d$. Indeed, there are
$q_1,\dots, q_n\in \mathbb Z$ such that
\begin{equation} \label{nod} \sum_{i=1}^nq_ia_i=1.\end{equation} After renumbering
of $\overline s_i$ we can assume that $q_i=-p_i<0$ for $i\leq k$ and
$q_i\geq 0$ for $i\geq k+1$. Denote by $M=a_1d\sum_{i=1}^k a_ip_i$
and for $j=0,1,\dots,a_1$ consider  elements
$$s_{0,j}=\prod_{i=1}^k\overline s_i^{(a_1-j)p_i}\cdot \prod_{i=k+1}^n\overline s_i^{jq_i}\in S_{O,{\bf{1}}}^G.$$
We have $$\tau(s_{0,j})=da_1\sum_{i=1}^kp_ia_i
+dj(-\sum_{i=1}^ka_ip_i+\sum_{i=k+1}^na_iq_i)=d(M +j)$$ for $0\leq
j\leq a_1$. Then $\tau(\overline s_1^m\cdot s_{0,j})=d(ma_1+M+j)$.
From this it is easy to see that $M$ satisfies the property that for
any $j\in \mathbb N$ there is an element $s\in S_{O,{\bf{1}}}^G$
with $\tau(s)=(M+j)d$, since $$\{ d(ma_1+M+j)\mid m\geq 0, 0\leq
j\leq a_1\}=d\mathbb N_{\geq M}.$$

To complete the proof of Proposition \ref{non-one}, note that
$N=M+t_0=M+\tau(s_1)$ is a desired constant. \qed \\

It is not difficult to give examples of groups $H$ satisfying
conditions of Proposition \ref{non-one}. For example, let
$H=SL_{p-1}(\mathbb Z_p)$ be the group of $(p-1)\times
(p-1)$-matrices with determinant $1$ over the finite field $\mathbb
Z_p$, $p\neq 2$. It is well-known that $H'=H$ and $Z(H)$, consisting
of scalar matrices, is a cyclic group of order $p-1$. For $i\neq j$
denote by $e_{i,j}$ the matrix whose entries are all zero except one
entry equal to one at the intersection of the $i$th row and $j$th
column. Put $t_{i,j}=e+e_{i,j}$, where $e$ is the identity matrix.
It is well known that the matrices $t_{i,j}$ (the transvections) are
all conjugate and that they generate the group $H=SL_{p-1}(\mathbb
Z_p)$. Therefore for  equipped group $(G,O)$, where
$G=PGL_{p-1}(\mathbb Z_p)$ and $O$ is the set of transvections,
almost all elements of $S_{O,{\bf{1}}}^G$ are not defined uniquely
by their type, that is, $S_{O,{\bf{1}}}^G$ (and, respectively,
$S_{O}$) is not a stable semigroup.

\section{Factorization semigroups over symmetric groups}
\subsection{Basic notations and definitions}
Let $\mathcal S_d$ be the symmetric group acting on the set $\{ 1,
\dots, d\}=I_d$. Remind that an element $\sigma =(i_1,\dots, i_k)\in
\mathcal S_d$ sending $i_1$ to $i_2$, $i_2$ to $i_3$, $\dots$,
$i_{k-1}$ to $i_k$, $i_k$ to $i_1$, and leaving fixed all over
elements of $I_d$ is called a {\it cyclic permutation} of {\it
length} $k$. A cyclic permutation of length $2$ is called a {\it
transposition}. Any cyclic permutation $\sigma=(i_{1},\dots,i_{k})$
is a product of $k-1$ transpositions:
\begin{equation}\label{cy2}\sigma= (i_k,i_{k-1})(i_{k-1},i_{k-2})\dots (i_2,i_1).\end{equation}
A factorization (\ref{cy2}) of $\sigma=(i_1,\dots ,i_k)$ is called
{\it canonical} if $i_1=\min_{1\leq j\leq k}\, i_j$.

As is well-known, any permutation $\sigma\in \mathcal S_d$,
$\sigma\neq \bf{1}$, can be represented as a product of cyclic
permutations:
\begin{equation} \label{cy1} \sigma =(i_{1,1},\dots,i_{k_1,1})(i_{1,2},\dots,i_{k_2,2})\dots
(i_{1,m},\dots,i_{k_m,m}),\end{equation} where $k_1\geq k_2\geq
\dots\geq k_m\geq 2$ and any two sets $\{
i_{1,j_1},\dots,i_{k_{j_1},j_1}\}$ and $\{
i_{1,j_2},\dots,i_{k_{j_2},j_2}\}$ have empty intersection if
$j_1\neq j_2$. If $\sigma$ is written in the form (\ref{cy1}), then
the ordered collection $t(\sigma)=[k_1,\dots, k_m]$ is called the
{\it type} of $\sigma$ and the number $l_t(\sigma)=\sum_{i=1}^m
k_i-m$ is called the {\it transposition length} of $\sigma$.

Note that for any $k_1\geq k_2\geq \dots\geq k_m\geq 2$ such that
$\sum k_j\leq d$ there is a permutation $\sigma$ of the type
$[k_1,\dots, k_m]$, and as is known, two permutations $\sigma_1$ and
$\sigma_2$ are conjugated in $\mathcal S_d$ if and only if
$t(\sigma_1)=t(\sigma_2)$. For a fixed type $t(\sigma)=[k_1,\dots,
k_m]$ a permutation
$$(1,\dots,k_1)(k_1+1,\dots,k_1+k_2)\dots(\sum_{i=1}^{m-1}k_i+1,\dots,\sum_{i=1}^{m}k_i)$$
is called the {\it canonical representative} of the type
$t(\sigma)$. The type $t(\sigma_1)=[k_{1,1},\dots, k_{m_1,1}]$ is
said to be {\it greater} than the type  $t(\sigma_2)=[k_{1,2},\dots,
k_{m_2,2}]$ if there is $l\geq 0$ such that $k_{1,i}=k_{2,i}$ for
$i\leq l$ and $k_{1,l+1}>k_{2,l+1}$ (here $k_{j,i}=0$ if $i>m_j$).
We say that a cyclic permutation $\sigma_1=({j_1},\dots, {j_{k_1}})$
is {\it greater} than a cyclic permutation $\sigma_2=({l_1},\dots,
{l_{k_2}})$ if either $t(\sigma_1)>t(\sigma_2)$ or if
$t(\sigma_1)=t(\sigma_2)$ then there is $r<k_1=k_2$ such that
${j_1}={l_1}, \dots, {j_r}={l_r}$, and ${j_{r+1}}>{l_{r+1}}$ in the
canonical factorizations of $\sigma_1$ and $\sigma_2$. Finally, we
say that a permutation $\sigma_1$ is {\it greater} than a
permutation $\sigma_2$ if either $t(\sigma_1)>t(\sigma_2)$ or if
$t(\sigma_1)=t(\sigma_2)$ and $\sigma_i=\sigma_{i,1}\dots
\sigma_{i,m}$, $i=1,2$, are cyclic factorizations, then there is $l$
such that $\sigma_{1,j}=\sigma_{2,j}$ for $j<l$ and
$\sigma_{1,l}>\sigma_{2,l}$. Denote by $\mathcal T =\{ t_1<t_2<\dots
<t_N\}$ the set of all types of permutations $\sigma\in \mathcal
S_d$.

By definition, the factorization semigroup  $\Sigma_d=S(\mathcal
S_d,\mathcal S_d)$ over the symmetric group $\mathcal S_d$ is
generated by the alphabet $X=\{ x_{\sigma} \mid \sigma \in \mathcal
S_d\}$. Let $s=x_{\sigma_1}\cdot\, .\, .\, .\, \cdot x_{\sigma_n}$
be an element of $\Sigma _d$. Applying relations (\ref{rel1}) and
(\ref{rel2}), we can assume that $t(\sigma_1)\leq \dots \leq
t(\sigma_n)$, then the sum $\tau(s)=\sum_{i=1}^N a_it_i$ is the {\it
type} of $s$, where $a_i$ is the number of factors $x_{\sigma_j}$,
$t(\sigma_j)=t_i$, interning in $s$.

For a subgroup $\Gamma$ of $\mathcal S_d$ denote
$\Sigma_d^{\Gamma}=\{ s\in \Sigma_d\mid \alpha(s)\in \Gamma\}$.

\subsection{Decompositions into products of transpositions}
Denote by $T_d$ the set of transpositions in $\mathcal S_d$. The
subsemigroup $S_{T_d}$ of $\Sigma_d$ is generated by $x_{(i,j)}$,
$1\leq i,j\leq d$, $i\neq j$, being subject to the relations
\begin{equation} \label{trrel}
\begin{array}{l}
x_{(i,j)}=x_{(j,i)} %&
\, \, \, \text{for all}\,\, \{ i,j\}_{ord} \subset I_d;
\\
x_{(i_1,i_2)}\cdot x_{(i_1,i_3)}=x_{(i_2,i_3)}\cdot
x_{(i_1,i_2)}=x_{(i_1,i_3)}\cdot x_{(i_2,i_3)} %&
\, \, \, \text{for all}\,\, \{i_1,i_2,i_3\}_{ord}\subset I_d;\\
x_{(i_1,i_2)}\cdot x_{(i_3,i_4)}=x_{(i_3,i_4)}\cdot x_{(i_1,i_2)} %&
\, \, \, \text{for all}\,\, \{ i_1,i_2,i_3,i_4\}_{ord} \subset I_d
\end{array}\end{equation}
(here $\{i_1,\dots , i_k\}_{ord}$ means a subset of $I_d$ consisting
of $k$ ordered elements, so that for any subset $\{ i_1,\dots,
i_k\}$ of $I_d$ we have $k!$ ordered subsets $\{ \sigma(i_1),\dots ,
\sigma(i_k)\}_{ord}$, $\sigma\in \mathcal S_k$).

Denote by $S_{T_d,\bf{1}}=S_{T_d}\cap \Sigma_{d,\bf{1}}$. By
Proposition \ref{simple} (4), the semigroup $\Sigma_{d,\bold 1}$ is
a subsemigroup of the center of $\Sigma_d$. In particular it is a
commutative semigroup.

It is easy to see that for each $\{ i,j\}\subset I_d$ the element
$s_{(i,j)}=x_{i,j}\cdot x_{i,j}=x_{(i,j)}^2$ belongs to
$S_{T_d,\bf{1}}$. The element
$$h_{d,g}=s_{(1,2)}^{g+1}\cdot s_{(2,3)}\cdot
.\, .\, .\cdot s_{(d-1,d)}\in S_{T_d,\bf{1}}\subset \Sigma_{d}$$ is
called a {\it Hurwitz element of genus} $g$.

\begin{lem} \label{chain1} For any ordered subset $\{ j_1,\dots , j_{k+1}\}_{ord}\subset I_d$
and for any $i$, $1\leq i\leq k$, the element $s=x_{(j_1,j_2)}\cdot
x_{(j_2,j_3)}\cdot \, .\, . \, .\cdot x_{(j_{k-1},j_k)}\cdot
x_{(j_i,j_{k+1})}\in S_{T_d}$ is equal to $$s_i= x_{(j_1,j_2)}\cdot
\, .\, .\, . \cdot x_{(j_{i-1},j_i)}\cdot x_{(j_i,j_{k+1})} \cdot
x_{(j_{k+1},j_{i+1})}\cdot x_{(j_{i+1},j_{i+2})} \cdot \, .\, . \,
.\cdot x_{(j_{k-1},j_k)}.$$
\end{lem}

\proof By (\ref{trrel}), we have (in each step of transformations
the underlining means that we will transform the underlined factors
and the result of transformation is written in brackets)
$$\begin{array}{l} s= x_{(j_1,j_2)}\cdot x_{(j_2,j_3)}\cdot \,
.\, . \, .\cdot \underline{ x_{(j_{i+1},j_{i+2})} \cdot \, .\, .\, . \cdot  x_{(j_{k-1},j_k)}\cdot x_{(j_i,j_{k+1})}}= \\
x_{(j_1,j_2)}\cdot \, .\, .\, . \cdot \underline{x_{(j_i,j_{i+1})}
\cdot (x_{(j_i,j_{k+1})}}\cdot x_{(j_{i+1},j_{i+2})} \cdot \, .\, .
\, .\cdot x_{(j_{k-1},j_k)})=
\\ x_{(j_1,j_2)}\cdot \, .\, .\, . \cdot x_{(j_{i-1},j_i)}(\underline{\cdot x_{(j_{i+1},j_{k+1})} \cdot
x_{(j_i,j_{i+1})}}) \cdot \, .\, . \, .\cdot x_{(j_{k-1},j_k)}=
\\ x_{(j_1,j_2)}\cdot \, .\, .\, . \cdot x_{(j_{i-1},j_i)}\cdot (x_{(j_i,j_{k+1})} \cdot
x_{(j_{k+1},j_{i+1})})\cdot x_{(j_{i+1},j_{i+2})} \cdot \, .\, . \,
.\cdot x_{(j_{k-1},j_k)}.
\end{array}$$
\qed

\begin{lem} \label{chain2} For any ordered subset $\{ j_1,\dots , j_{k}\}_{ord}\subset I_d$
and for any $i$, $1\leq i\leq k$, the element $s=x_{(j_1,j_2)}\cdot
x_{(j_2,j_3)}\cdot \, .\, . \, .\cdot x_{(j_{k-1},j_k)}\cdot
x_{(j_i,j_k)}\in S_{T_d}$, where $k\leq d-1$, is equal to
$s_i=x_{(j_1,j_2)}\cdot\, .\, .\, .\, \cdot x_{(j_{i-1},j_i)}\cdot
x_{(j_{i+1},j_{i+2})}\cdot\, .\, .\, .\, \cdot
x_{(j_{k-1},j_k)}\cdot x_{(j_i,j_{i+1})}^2$.
\end{lem}

\proof By (\ref{trrel}), we have
$$\begin{array}{l} s=x_{(j_1,j_2)}\cdot x_{(j_2,j_3)}\cdot \,
.\, . \, .\cdot \underline{x_{(j_{k-1},j_k)}\cdot x_{(j_i,j_k)}}= \\
x_{(j_1,j_2)}\cdot x_{(j_2,j_3)}\cdot \,
.\, . \, .\cdot \underline{x_{(j_{k-2},j_{k-1})}\cdot (x_{(j_i,j_{k-1})}}\cdot x_{(j_{k-1},j_k)})= \dots =\\
x_{(j_1,j_2)}\cdot \, .\, .\, . \cdot x_{(j_{i-1},j_i)} \cdot
x_{(j_i,j_{i+1})}\cdot (x_{(j_i,j_{i+1})}\cdot
x_{(j_{i+1},j_{i+2})}) \cdot \, .\, . \, .\cdot x_{(j_{k-1},j_k)}=
\\ x_{(j_1,j_2)}\cdot \, .\, .\, . \cdot x_{(j_{i-1},j_i)} \cdot
\underline{x_{(j_i,j_{i+1})}^2\cdot x_{(j_{i+1},j_{i+2})}} \cdot \,
.\, . \, .\cdot x_{(j_{k-1},j_k)}= \\ x_{(j_1,j_2)}\cdot \, .\, .\,
. \cdot x_{(j_{i-1},j_i)}\cdot (x_{(j_{i+1},j_{i+2})} \cdot
\underline{x_{(j_i,j_{i+1})}^2) \cdot
x_{(j_{i+2},j_{i+3})}}\cdot \, .\, . \, .\cdot x_{(j_{k-1},j_k)}=\dots = \\
x_{(j_1,j_2)}\cdot\, .\, .\, .\, \cdot x_{(j_{i-1},j_i)}\cdot
x_{(j_{i+1},j_{i+2})}\cdot\, .\, .\, .\, \cdot
(x_{(j_{k-1},j_k)}\cdot x_{(j_i,j_{i+1})}^2)=s_i.
\end{array}$$
\qed
\begin{lem} \label{chain4} The following equalities:
\begin{equation} \label{comm1} x_{(i_1,i_2)}^2\cdot x_{(i_2,i_3)}= x_{(i_2,i_3)}\cdot
x_{(i_1,i_3)}^2 = x_{(i_1,i_3)}^2\cdot x_{(i_2,i_3)}=
x_{(i_2,i_3)}\cdot x_{(i_1,i_2)}^2;\end{equation}
\begin{equation}\label{retd1}
x_{(i_1,i_2)}^2\cdot x_{(i_2,i_3)}^2=x_{(i_1,i_2)}^2\cdot
x_{(i_1,i_3)}^2=x^2_{(i_2,i_3)}\cdot x^2_{(i_1,i_3)} \end{equation}
hold for all ordered triples $\{ i_1,i_2,i_3\}_{ord}\subset I_d$;
and
\begin{equation}\label{retd2}
x^2_{(i_1,i_2)}\cdot x^2_{(i_3,i_4)}=x^2_{(i_3,i_4)}\cdot
x^2_{(i_1,i_2)}
\end{equation} hold for all ordered $4$-tuples $\{ i_1,i_2,i_3,i_4\}_{ord}\subset
I_d$.
\end{lem}

\proof We will check only two of three equalities (\ref{comm1}),
since the inspection of the other equalities is similar. By
(\ref{trrel}), we have
$$\begin{array}{l} x_{(i_1,i_2)}^2\cdot x_{(i_2,i_3)}=x_{(i_1,i_2)}\cdot \underline{x_{(i_1,i_2)}\cdot
x_{(i_2,i_3)}}=  \underline{x_{(i_1,i_2)}\cdot (x_{(i_2,i_3)}}\cdot
x_{(i_1,i_3)})=
\\ (x_{(i_2,i_3)}\cdot x_{(i_1,i_3)})\cdot
x_{(i_1,i_3)}=x_{(i_2,i_3)}\cdot x_{(i_1,i_3)}^2. \end{array} $$

Similarly,
$$\begin{array}{l} x_{(i_1,i_2)}^2\cdot x_{(i_2,i_3)}=x_{(i_1,i_2)}\cdot \underline{x_{(i_1,i_2)}\cdot
x_{(i_2,i_3)}}=  \underline{x_{(i_1,i_2)}\cdot (x_{(i_1,i_3)}}\cdot
x_{(i_1,i_2)})=
\\ (x_{(i_2,i_3)}\cdot x_{(i_1,i_2)})\cdot
x_{(i_1,i_2)}=x_{(i_2,i_3)}\cdot x_{(i_1,i_2)}^2. \end{array} $$
\qed \\

The following lemma is a particular case of Lemma \ref{fried}.

\begin{lem} \label{chain3} For any ordered subset $\{ j_1,\dots , j_{k}\}_{ord}\subset I_d$
the following equality:
$$x_{(j_1,j_2)}^2\cdot x_{(j_1,j_2)}\cdot x_{(j_2,j_3)}\cdot \, .\, .
\, .\cdot x_{(j_{k-1},j_k)}= x_{(j_i,j_l)}^2\cdot x_{(j_1,j_2)}\cdot
x_{(j_2,j_3)}\cdot \, .\, . \, .\cdot x_{(j_{k-1},j_k)}$$ holds,
where $1\leq i<l\leq k$.
\end{lem}

To each word $w(\overline{x_{(i,j)}})=x_{(i_1,j_1)}\dots
x_{(i_m,j_m)}\in W=W(T_d)$, let us associate a graph
$\widetilde{\Gamma}_w$ consisting of $d$ vertices $v_i$, $1\leq
i\leq d$, the set of edges is in one to one correspondence with the
collection of letters incoming in $w$ so that two vertices $v_i$ and
$v_j$ are connected by an edge if the letter $x_{(i,j)}$ is
contained in $w$, in particular the number of edges connecting
vertices $v_i$ and $v_j$ coincides with the number of entry of the
letter $x_{(i,j)}$ in $w$. The edges of the graph
$\widetilde{\Gamma}_w$ are numbered according to the position of the
corresponding letter in $w$. Denote by $V_{\text{iso}}$ the set of
isolated vertices of $\widetilde{\Gamma}_w$ (that is, a vertex $v_i$
is {\it isolated} if it is not connected by an edge with some other
vertex of $\widetilde{\Gamma}_w$) and put
$\Gamma_w=\widetilde{\Gamma}_w\setminus V_{\text{iso}}$.
\begin{lem} \label{graf} For any $s\in S_{T_d}$ and for any $w_1,w_2\in W(s)$ the graphs
$\Gamma_{w_1}$ and $\Gamma_{w_2}$ have the same sets of vertices
$V(s)=V(\Gamma_{w_1})=V(\Gamma_{w_2})$.
\end{lem}

\proof It is easily follows from relations (\ref{trrel}). \qed
\begin{prop} \label{canonic}  Let $s\in S_{T_d}$ be of length $k\leq d-1$.
Then  $\alpha(s)\in \mathcal S_d$ is a cyclic permutation of length
$k$ if and only if $s$ satisfies the following condition:

$$\text{there is a word}\, \,  w\in W(s)\, \, \text{ those graph}\, \,  \Gamma_{w}\, \, \text{is a tree}.
\qquad \qquad (*)$$ Moreover, an element $s$ satisfying condition
$(*)$ is uniquely defined by the cyclic permutation $\alpha(s)$.
\end{prop}
\proof Let us show that if $s$ satisfies condition $(*)$, then there
are exactly $k=ln(s)$ words $w_1,\dots, w_k\in W(s)$ such that
$\Gamma_{w_i}$ are simple paths if we go along the edges according
to their numbering. Indeed, it is easy to see that Lemma
\ref{chain1} implies the existence of a word
$w_1=x_{(i_1,i_2)}x_{(i_2,i_3)}\dots x_{(i_{k-1},i_k)}$ whose graph
$\Gamma_{w_1}$ is a simple path. Let us show that if we move the
letter $x_{(i_{k-1},i_k)}$ to the left then we again obtain a word
$w_2$ defining the same element $s$ and such that $\Gamma_{w_2}$ is
a simple path. Indeed, we have
$$\begin{array}{ll}
s= & x_{(i_1,i_2)}\cdot \, .\, .\, .\, \cdot
\underline{x_{(i_{k-2},i_{k-1})}\cdot
x_{(i_{k-1},i_k)}}=  \\
& x_{(i_1,i_2)}\cdot \, .\, .\, .\, \cdot
\underline{x_{(i_{k-3},i_{k-2})}\cdot (x_{(i_{k-2},i_{k})}}\cdot
x_{(i_{k-2},i_{k-1})})= \dots =  \\
& (x_{(i_1,i_k)}\cdot x_{(i_1,i_2)})\cdot \, .\, .\, .\, \cdot
x_{(i_{k-2},i_{k-1})}.\end{array}$$ Repeating such transformations
$k$ times, we find desired words $w_1,\dots , w_k$.

We have $\alpha(s) =(i_1,i_2)\dots (i_{k-2},i_{k-1})(i_{k-1},i_k)$
is a cyclic permutation of length $k$. On the other hand, if
$\sigma\in \mathcal S_d$ is a cyclic permutation of length $k$ then
it can be represented as a product of $k-1$ transpositions
$\sigma=(i_1,i_2)\dots (i_{k-2},i_{k-1})(i_{k-1},i_k)$ and,
obviously, that $\alpha(s)=\sigma$ for $s= x_{(i_1,i_2)}\cdot \, .\,
.\, .\, \cdot x_{(i_{k-2},i_{k-1})}\cdot x_{(i_{k-1},i_k)}$ and the
graph $\Gamma_{x_{(i_1,i_2)}\dots x_{(i_{k-2},i_{k-1})}
x_{(i_{k-1},i_k)}}$ satisfies condition $(*)$.

Now if we fix a set $\{ i_1,\dots, i_k\}\subset I_d$ then there are
exactly $(k-1)!$ distinct cyclic permutations in $\mathcal S_d$ of
length $k$ cyclicly permuting the elements of the set $\{ i_1,\dots,
i_k\}$. On the other hand, there are exactly $k!$ distinct simple
paths connecting the vertices $v_{i_1},\dots ,v_{i_k}$. Therefore,
the elements $s$ satisfying condition $(*)$ are defined uniquely by
the cyclic permutations $\alpha(s)$. \qed

\begin{thm} \label{fuc} For any $s\in S_{T_d}$ the difference
$ln(s)-l_t(\alpha(s))$ is a non-negative even number and there are
elements $\widetilde s\in S_{T_d}$ and $\overline s\in
S_{T_d,\bf{1}}$ such that $s=\widetilde s\cdot \overline s$, the
length $ln(\widetilde s)=l_t(\alpha(s))$ and $\alpha(\widetilde
s)=\alpha(s)$.

If $s\in S_{T_d}^{\mathcal S_d}$ and $ln(s)\geq
l_t(\alpha(s))+2(d-1)$, then one can find a factorization
$s=\widetilde s\cdot \overline s$, where $\overline s=h_{d,g}$ with
$g=\frac{1}{2}(ln(s)- l_t(\alpha(s)))-d+1$ and $\widetilde s$ is
such that $ln(\widetilde s)=l_t(\alpha(s))$, $\alpha(\widetilde
s)=\alpha(s)$, moreover, $\widetilde s$ is defined uniquely by
$\alpha(s)$.
\end{thm}
\proof  Consider the graph $\Gamma_w$ of some $w\in W(s)$. It splits
into the disjoint union of its connected components:
$\Gamma_w=\Gamma_{w,1}\sqcup\dots\sqcup \Gamma_{w,l}$. It is easily
follows from (\ref{trrel}) that
$s=\varphi(w_1(\overline{x_{(i,j)}}))\cdot \, .\, .\, .\, \cdot
\varphi(w_l(\overline{x_{(i,j)}}))$, where
$w_i(\overline{x_{(i,j)}})$ is a word in letters $x_{(i,j)}$'s such
that $\Gamma_{w_i}=\Gamma_{w,i}$. Let $s_i=\varphi(w_i)\in S_{T_d}$
be an element defined by the word $w_i$. We have $(\mathcal
S_d)_{s_i}\cap (\mathcal S_d)_{s_j}=\bf 1$ for $i\neq j$, in
particular, $s_i\cdot s_j=s_j\cdot s_i$.  Applying Lemma
\ref{chain1}, it is easy to see that for each $i$ we can find a
representation of $s_i$ as a word in letters $x_{(i,j)}$'s such that
$$s_i=x_{(j_{1,i},j_{2,i})}\cdot \,
.\, . \, .\cdot x_{(j_{k_i-1,i},j_{k_i,i})}\cdot s_{i,1}$$ and the
set $\{ v_{j_1,i},\dots, v_{j_{k_i},i}\}$ is the complete set of the
vertices of $\Gamma_{w_i}$.

Let $x_{(j_a,j_b)}$, $a<b$, be the first factor of $s_{i,1}$  if
$s_{i,1}\neq x_{\bf{1}}$. Then it follows from relations
(\ref{trrel}) and Lemma \ref{chain2} that $s_i$ can be written in
the form : $s_i=s_i'\cdot x_{(j_a,j_b)}^2$. Note that
$x_{(j_a,j_b)}^2\in S_{T_d,\bf{1}}$ and
$ln(s_i')=ln(s_i)-2<ln(s_i)$, that is, we obtain that $s$ can be
written in the form: $s=\widetilde s_1\cdot \overline s_1$, where
$ln(\widetilde s_1)<ln(s)$ and $\overline s_1\in S_{T_d,\bf{1}}$, in
addition, $\alpha(\widetilde s_1)=\alpha(s)$, since $\overline
s_1\in S_{T_d,\bf{1}}$. Repeating, if necessary,  these arguments
for $\widetilde s_1$, $\dots$, as a result we obtain that $s$ can be
written in the form: $s=\widetilde s\cdot \overline s$, where
$\overline s\in S_{T_d,\bf{1}}$ is a product of some squares of
$x_{(i,j)}$'s and $\widetilde s=s_1\cdot \, .\, .\, . s_m\in
S_{T_d}$, where for $1\leq i\leq m$ the elements
$s_i=x_{(j_{1,i},j_{2,i})}\cdot \, .\, .\, . \cdot
x_{(i_{k_i-1,i},j_{k_i,i})}$ are such that the subsets $\{
j_{1,i},\dots ,j_{k_i,i}\}$ and $\{ j_{1,l},\dots ,j_{k_l,l}\}$ of
$I_d$ have the empty intersection for $i\neq l$. Therefore
$$\alpha(s)=\alpha(\widetilde s)=(j_{k_1,1},\dots,j_{1,1})\dots
(j_{k_m,m},\dots,j_{1,m})$$ and hence $ln(\widetilde
s)=l_t(\alpha(s))$.

Therefore, by Proposition \ref{canonic},  the elements $s_i$ are
defined uniquely (up to renumbering) by $\alpha(s_i)$.

Now let $s=\widetilde s\cdot \overline s\in S_{T_d}^{\mathcal S_d}$
with $ln(s)\geq l_t(\alpha(s))+2(d-1)$, where $\overline s\in
S_{T_d\bf{1}}$ is a product of some squares of $x_{(i,j)}$'s and
$\widetilde s$ is such that $$\alpha(s)=\alpha(\widetilde
s)=(j_{1,1},\dots,j_{k_1,1})\dots (j_{1,m},\dots,j_{k_m,m})$$ and
$ln(\widetilde s)=l_t(\alpha(s))$. Note that $ln(\overline s)\geq
2(d-1)$, since $ln(\widetilde s)=l_t(\alpha(s))$.

Consider the graphs $\Gamma_{\widetilde w}$, $\Gamma_{\overline w}$,
and $\Gamma_{\widetilde w\overline w}$, where $\widetilde w\in
W(\widetilde s)$, $\overline w\in W(\overline s)$, and $\widetilde
w\overline w\in W(s)$. Let us show that there is a factorization of
$s=\widetilde s\cdot \overline s$ such that $V_{\overline s}=I_d$.
First of all, we have $V_s=I_d$, since $(\mathcal S_d)_s=\mathcal
S_d$. Assume that $V_{\overline s}\neq I_d$ for some factorization
of $s=\widetilde s\cdot \overline s$ and let $\overline
s=\varphi(\overline w(\overline{x_{(i,j)}^2}))$ and $\widetilde
s=\varphi(\widetilde w(\overline{x_{(i,j)}}))$. Since $ln(\overline
s)\geq 2(d-1)$, it follows from Lemma \ref{chain4} that there is a
connected component $\Gamma_1$ of $\Gamma_{\overline w}$ such that
for any pair of vertices $v_{i_1},v_{i_2}\in \Gamma_1$ we can find a
word $\overline w\in W(\overline s)$ such that $\overline s=
(x_{(i_1,i_2)}^2)^2\cdot \overline s'$. Next, since $V_s=I_d$, then
there is a pair $v_{i_0},v_{i_2}\in V_{\widetilde s}$ such that
$v_{i_0}\not\in V_{\overline s}$, $v_{i_2}\in V_{\overline s}$, and
$\widetilde s=\widetilde s'\cdot x_{(i_0,i_2)}$. By Lemma
\ref{chain4}, we have
$$s=\widetilde s\cdot \overline s=\widetilde s'\cdot
x_{(i_0,i_2)}\cdot x_{(i_1,i_2)}^2\cdot x_{(i_1,i_2)}^2\cdot
\overline s'=\widetilde s'\cdot x_{(i_0,i_2)}\cdot
x_{(i_0,i_1)}^2\cdot x_{(i_1,i_2)}^2\cdot \overline s'=\widetilde
s\cdot \overline s_1, $$ where either $V_{\overline
s_1}=V_{\overline s}\cup \{ i_0\}$ or for a word $\overline w_1\in
W(\overline s_1)$ the number of connected components of the graph
$\Gamma_{\overline w_1}$ is strictly less than the number of
connected components of $\Gamma_{\overline w}$. Repeating this
transformation several times, as a result we obtain a factorization
$s=\widetilde s\cdot \overline s$, such that $V_{\overline s}=I_d$.
Now, to complete the proof of Theorem \ref{fuc} it suffices to apply
once more Lemma \ref{chain4}. \qed
\begin{prop} \label{redegen} There is a unique homomorphism $r:\Sigma_d\to
S_{T_d}$ such that
\begin{itemize}
\item[$(i)$] $\alpha (r(x_{\sigma}))=\sigma$ for $\sigma\in \mathcal S_d$,
\item[$(ii)$] $ln(r(x_{\sigma}))=l_t(\sigma)$,
\item[$(iii)$] $r_{|S_{T_d}}=Id$.
\end{itemize} \end{prop}

\proof Each element $\sigma\in \mathcal S_d$, $\sigma\neq \bf{1}$,
can be factorized into a product of pairwise commuting cycles:
$\sigma=\sigma_1\dots \sigma_m$ and such a factorization is unique
up to permutations of factors. According to Proposition
\ref{canonic}, each of these cyclic permutations $\sigma_i$ defines
uniquely an element $s_i\in S_{T_d}$ such that $ln(s_i)=k_i-1$ and
$\alpha(s_i)=\sigma_i$, where $k_i$ is the length of the cycle
$\sigma_i$, and therefore the product $s(\sigma)=s_1\cdot \, .\, .\,
.\, \cdot s_m\in S_{T_d}$ is defined uniquely by $\sigma$. It is
easy to see that the map $\sigma \mapsto s(\sigma)$ defines a
homomorphism $r:\Sigma_d\to S_{T_d}$ given by
$r(x_{\sigma})=s(\sigma)$ on the set of generators of $\Sigma_d$. It
is obvious that
$ln_t(s)=ln(r(s))$ and $r_{|S_{T_d}}=Id$. \qed \\

The homomorphism $r:\Sigma_d\to S_{T_d}$ defined in Proposition
\ref{redegen} is called the {\it regenerating} homomorphism and the
number $ln_t(s)=ln(r(s))$ is called the {\it transposition length}
of $s\in \Sigma_d$.

\subsection{Decompositions of the unity into products of transpositions}
Let us consider the semigroup $S_{T_d,\bf{1}}$.
\begin{thm} \label{std1} The semigroup $S_{T_d,\bf{1}}$ is
 commutative and it is generated by the elements $s_{(i,j)}=x_{(i,j)}^2$, $\{ i,j\}\subset I_d$,
being subject to the relations
\begin{equation}\label{retd11}
s_{(i_1,i_2)}\cdot s_{(i_2,i_3)}=s_{(i_1,i_2)}\cdot
s_{(i_1,i_3)}=s_{(i_2,i_3)}\cdot s_{(i_1,i_3)} \end{equation} for
all ordered triples $\{ i_1,i_2,i_3\}_{ord}\subset I_d$ and
\begin{equation}\label{retd21}
s_{(i_1,i_2)}\cdot s_{(i_3,i_4)}=s_{(i_3,i_4)}\cdot s_{(i_1,i_2)}
\end{equation} for all ordered $4$-tuples $\{ i_1,i_2,i_3,i_4\}_{ord}\subset
I_d$. Moreover, any element $s\in S_{T_d,\bf{1}}$ has a normal form,
that is, it can be uniquely written in the form
$$s=(s^{k_1}_{(i_{1,1},i_{2,1})}\cdot s_{(i_{2,1},i_{3,1})}\cdot \, .\,
.\, .\, \cdot  s_{(i_{j_1-1,1},i_{j_1,1})})\cdot \, .\, .\, .\,
\cdot (s^{k_n}_{(i_{1,n},i_{2,n})}\cdot s_{(i_{2,n},i_{3,n})}\cdot
\, .\, .\, .\, \cdot  s_{(i_{j_n-1,n},i_{j_n,n})}),$$ where $1\leq
i_{1,1}<i_{1,2}<\dots <i_{1,n}\leq d-1$,\, \, $k_l\in\mathbb N$ for
$l=1,\dots, n$, the sets $M_l=\{ i_{1,l}<i_{2,l},\dots
<i_{j_l,l}\}$, $1\leq l\leq n$, are subsets of $I_d$ of cardinality
$j_l\geq 2$ such that $M_{l_1}\cap M_{l_2}=\emptyset$ for $l_1\neq
l_2$.
\end{thm}

\proof It follows from Theorem \ref{fuc} that $S_{T_d,\bf{1}}$ is
generated by $s_{(i,j)}$'s. By Lemma \ref{chain4},  the elements
$s_{(i,j)}$ satisfy relations (\ref{retd11}) and (\ref{retd21}).

Like in the proof of Theorem \ref{fuc}, for each
$s=s_{(j_1,j_2)}\cdot \, .\, .\, .\, \cdot s_{(j_{m-1},j_m)}$ we can
associate a graph $\Gamma_w$, where $w$ is a word in letters
$s_{(i,j)}$ representing the element $s$. The graph $\Gamma_w$
splits into the disjoint union of its connected components:
$\Gamma_w=\Gamma_{w,1}\sqcup\dots\sqcup \Gamma_{w,n}$. It is easily
follows from (\ref{trrel}) that $w=w_1(\overline{s_{(i,j)}})\dots
w_n(\overline{s_{(i,j)}})$, where $w_l(\overline{s_{(i,j)}})$ is a
word in letters $s_{(i,j)}$'s such that $\Gamma_{w_l}=\Gamma_{w,l}$.
Let $s_l\in S_{T_d,\bf{1}}$ be an element defined by the word $w_l$,
that is, $s_l=\varphi(w_l)$.

It is easily  follows from relations (\ref{retd11}) and
(\ref{retd21}) that each element $s_l$ can be uniquely written in
the form
\begin{equation} \label{xxx} s_l=s^{k_l}_{(i_{1,l},i_{2,l})}\cdot s_{(i_{2,l},i_{3,l})}\cdot
\, .\, .\, .\, \cdot  s_{(i_{j_l-1,l},i_{j_l,l})},\end{equation}
where the set $M_l=\{ i_{1,l}<i_{2,l},\dots <i_{j_l,l}\}$, $1\leq
l\leq n$, is in one to one correspondence with the set of vertices
of the connected component $\Gamma_{w,l}$ of the graph $\Gamma_w$.
\qed
\begin{rem} {\rm Note that the element $s^{k_l}_{(i_{1,l},i_{2,l})}\cdot s_{(i_{2,l},i_{3,l})}\cdot
\, .\, .\, .\, \cdot  s_{(i_{j_l-1,l},i_{j_l,l})}$ in (\ref{xxx}) is
the Hurwitz element $h_{j_l,k_l-1}$ of the semigroup
$S_{T_{j_l},\bf{1}}$ if we consider $S_{T_{j_l},\bf{1}}$ as a
subsemigroup of} $S_{T_{d,\bf{1}}}$ {\rm and the embedding is
defined by the natural embedding} $M_l\hookrightarrow I_d$.
\end{rem}

\begin{prop} \label{hur}
The Hurwitz element $h_{d,g}$ belongs to the center of the semigroup
$\Sigma_d$ and it is fixed under the conjugation action of $\mathcal
S_m$ on $\Sigma_d$.

For $h_{d,g_1},\, h_{d,g_2}$ we have  $$h_{d,g_1}\cdot
h_{d,g_2}=h_{d,g_1+g_2+d-1}.$$
\end{prop}

\proof The first part of Proposition follows from Proposition
\ref{simple}, since, on the one hand, $\alpha (h_{d,g})=\bold 1$ and
the transpositions $(i,i+1)$, $i=1,\dots , d-1$, generate the group
$(\mathcal S_d)_{h_{d,g}}$. On the other hand, they generate the
symmetric group $\mathcal S_d$.

The second part of Proposition follows from Theorem \ref{fuc}.
 \qed \\

Moreover, as a corollary of Theorems \ref{fuc} and \ref{std1} we
obtain that a Hurwitz element $h_{d,g}$ is defined uniquely in the
semigroup $S_{T_d}$ by its length and  the following two conditions.
\begin{cor} {\bf (Clebsch -- Hurwitz Theorem)} {\rm (\cite{Cl})} Let an element $s\in S_{T_d}$
satisfy the following conditions \begin{itemize}
\item[($i$)] $(\mathcal S_d)_s=\mathcal S_d$;
%\item[($ii$)] $ln(s)\geq 2(d-1)$;
\item[($ii$)] $\alpha(s)=\mathbf 1$.
\end{itemize} Then $ln(s)\geq 2(d-1)$ and $s=h_{d,g}$, where $g=\frac{ln(s)}{2}-d+1$.
\end{cor}

\subsection{Factorizations in the symmetric groups (general
case)} In this subsection we will prove the following generalization
of Theorem \ref{fuc}.
\begin{thm} \label{fuc1} Let $s=x_{\sigma_{1}}\cdot \, .\, .\,
.\, \cdot x_{\sigma_{m}}\cdot \overline s\in \mathcal S_d$, where
$\overline s\in S_{T_d}$. For $j=1,\dots , m$, denote by
$\sigma_{j,0}$ the canonical representative of the type
$t(\sigma_j)$ and by
$$\sigma=\sigma(s)=(\sigma_{1,0}\dots \sigma_{m,0})^{-1}\alpha(s).$$
If $s\in \Sigma_d^{\mathcal S_d}$ and $ln(\overline s)=k\geq
3(d-1)$, then
$$s=x_{\sigma_{1,0}}\cdot \, .\, .\, .\, \cdot x_{\sigma_{m,0}}\cdot
r(x_{\sigma})\cdot h_{d,g},$$ where
$g=\frac{k-ln_t(x_{\sigma})}{2}-d+1$.
\end{thm}
\proof  Let us show that there is a factorization $$
s=x_{\sigma_1'}\cdot \, .\, .\, .\, \cdot x_{\sigma_m'}\cdot
x_{(i_1,j_1)}\cdot \, .\, .\, .\, \cdot
x_{(i_k,j_k)}=x_{\sigma_1'}\cdot \, .\, .\, .\, \cdot
x_{\sigma_m'}\cdot \overline s_1$$ such that
$t(\sigma_i)=t(\sigma_i')$ for $i=1,\dots, m$ and the set
$V_{\overline s_1}$ of vertices of the graph $\Gamma_{\overline
w_1}$ of the word $\overline w_1=x_{(i_1,j_1)}\dots x_{(i_k,j_k)}\in
W(\overline s_1)$ coincides with the set $I_d$.

Indeed, let $w\in W(\overline s)$ and assume that $V_{\overline
s}\neq I_d$. Since $ln(\overline s)\geq 3(d-1)$, then there is a
connected component $\Gamma_1$ of the graph $\Gamma_{w}$ such that
the number of its edges is greater than the number of its vertices.
Then it follows from the proof of Theorem \ref{fuc} that for any
$v_{i_1},v_{i_2}$ belonging to the set $V(\Gamma_1)$ of vertices of
$\Gamma_1$ there is a word $w'\in W$  such that $\overline s=
x_{(i_1,i_2)}^2\cdot \varphi(w')$ and the vertices of $V(\Gamma_1)$
belong to one and the same connected component of
$\Gamma_{x_{i_1,i_2}^2w'}$. Next, since $(\mathcal S_d)_{s}=\mathcal
S_d$, then there is $\sigma_{l}$ for some $l$, $1\leq l\leq m$, such
that $\sigma_{l}(i_1,i_2)\sigma_{l}^{-1}=(i_0,j_0)$, where either
$v_{i_0}$ or $v_{j_0}$ (but not both) does not belong $V(\Gamma_1)$.
Without loss of generality, we can assume that $l=m$. We have
$$\begin{array}{l} s=x_{\sigma_{1}}\cdot \, .\, .\, .\, \cdot
x_{\sigma_{m}}\cdot \overline s=  x_{\sigma_{1}}\cdot \, .\, .\, .\,
\cdot x_{\sigma_{m}}\cdot x_{(i_1,i_2)}^2\cdot \varphi(w')= \\
x_{\sigma_{1}}\cdot \, .\, .\, .\, \cdot x_{\sigma_{m-1}} \cdot
x_{(i_0,j_0)}\cdot
x_{\sigma_{m}}\cdot x_{(i_1,i_2)}\cdot \varphi(w')= \\
x_{\sigma_{1}}\cdot \, .\, .\, .\, \cdot x_{\sigma_{m-1}}  \cdot
\rho((i_0,j_0))(x_{\sigma_{m}})\cdot x_{(i_0,j_0)} \cdot
x_{(i_1,i_2)}\cdot \varphi(w')= \\ x_{\sigma_{1}}\cdot \, .\, .\,
.\, \cdot x_{\sigma_{m-1}}  \cdot \rho((i_0,j_0))(x_{\sigma_{m}})
\cdot \varphi(w''),
\end{array}$$
where $w''=x_{(i_0,j_0)}x_{(i_1,j_1)}w'$ such that either the set of
vertices of $\Gamma_{w''}$ strictly contains the set $V_{\overline
s}$ or the number of connected components of $\Gamma_{w''}$ is
strictly less than the one of $\Gamma_{w'}$.

Repeating such transformations several times, as a result we obtain
a factorization of $s$ of the form
$$s=x_{\sigma_1'}\cdot \, .\, .\, .\, \cdot
x_{\sigma_m'}\cdot \overline s_1$$ such that $\overline s_1\in
S_{T_d}$ and $V_{\overline s_1}=I_d$, and $t(\sigma_j')=t(\sigma_j)$
for $j=1,\dots, m$. For this factorization we have $(\mathcal
S_d)_{\overline s_1}=\mathcal S_d$ and $ln(\overline s_1)\geq
3(d-1)$.

To complete the proof of Theorem \ref{fuc1} we will use induction by
$m$. For $m=0$ Theorem \ref{fuc1} follows from Theorem \ref{fuc}.

Let $m=1$. By Theorem \ref{fuc} we have $\overline s_1 =h_{d,0}\cdot
\overline s'$ for some $\overline s'\in S_{T_d}$.
\begin{lem} \label{hhhh} For any disjoint union  $\{ i_{1,1},\dots, i_{k_1,1}\}
\sqcup \dots \sqcup \{ i_{1,n},\dots, i_{k_n,n}\}$ of ordered
subsets of $I_d$ the element $h_{d,0}$ can be represented as a
product
$$ h_{d,0}=(x_{(i_{1,1},i_{2,1})}\cdot \,. \, .\, .\, \cdot x_{(
i_{k_1-1,1},i_{k_1,1})})\cdot \, .\, .\, .\, \cdot
(x_{(i_{1,n},i_{2,n})}\cdot \,. \, .\, .\, \cdot x_{(
i_{k_n-1,n},i_{k_n,n})})\cdot \overline h,$$ where $\overline h$ is
an element of $S_{T_d}^{\mathcal S_{d}}$.
\end{lem}
\proof The subgroup $S_{T_{d,\bf{1}}}$ is commutative and the
element $h_{d,0}$ is invariant under the conjugation action of
$\mathcal S_d$, therefore $h_{d,0}$ can be written in the form
$$ h_{d,0}=(s_{(i_{1,1},i_{2,1})}\cdot \,. \, .\, .\, \cdot s_{(
i_{k_1-1,1},i_{k_1,1})})\cdot \, .\, .\, .\, \cdot
(s_{(i_{1,n},i_{2,n})}\cdot \,. \, .\, .\, \cdot s_{(
i_{k_n-1,n},i_{k_n,n})})\cdot \widetilde h,$$ where $\widetilde h$
is an element of $S_{T_{d,\bf{1}}}$. We have
$$\begin{array}{l} s_{(i_{1,j},i_{2,j})}\cdot \,. \, .\, .\, \cdot s_{(
i_{k_j-1,j},i_{k_j,j})}=x^2_{(i_{1,j},i_{2,j})}\cdot \,. \, .\, .\,
\cdot x^2_{( i_{k_j-1,j},i_{k_j,j})}= \\
x_{(i_{1,j},i_{2,j})}\cdot (x^2_{(i_{2,j},i_{3,j})}\cdot \,. \, .\,
.\, \cdot x^2_{( i_{k_j-1,j},i_{k_j,j})})\cdot
x_{(i_{1,j},i_{2,j})}=\dots = \\
(x_{(i_{1,j},i_{2,j})}\cdot \,. \, .\, .\, \cdot x_{(
i_{k_j-1,j},i_{k_j,j})}) \cdot (x_{(i_{k_j-1,j},i_{k_j,j})}\cdot \,.
\, .\, .\, \cdot x_{( i_{1,j},i_{2,j})})
\end{array}$$
and the elements $x_{(i_{l_1,j_1},i_{l_1+1,j_1})}$ and
$x_{(i_{l_2,j_2},i_{l_2+1,j_2})}$ commute if $j_1\neq j_2$. Now to
complete the proof of Lemma, note that $V_{s_j}=V_{\overline s_j}$,
where $s_j=s_{(i_{1,j},i_{2,j})}\cdot \,. \, .\, .\, \cdot s_{(
i_{k_j-1,j},i_{k_j,j})}$ and $\overline
s_j=x_{(i_{k_j-1,j},i_{k_j,j})}\cdot \,. \, .\, .\, \cdot x_{(
i_{1,j},i_{2,j})}$. Therefore $V_{\overline h}=I_d$ for
$\overline h =(\prod \overline s_i)\cdot \widetilde h$. \qed \\

For the canonical representative $\sigma_{m,0}$ of the type
$t(\sigma_{m})$ there is $\overline{\sigma}_m\in \mathcal S_d$ such
that
$\sigma_{m,0}=\overline{\sigma}_m^{-1}\sigma_{m}'\overline{\sigma}_m$.
The permutation $\overline{\sigma}_m$ can be factorized into the
product of cyclic permutations and each cyclic permutation can be
factorized into the product of transpositions:
$$\overline{\sigma}_m=((i_{1,1},i_{2,1})\dots (
i_{k_1-1,1},i_{k_1,1}))\dots ((i_{1,n},i_{2,n})\dots (
i_{k_n-1,n},i_{k_n,n})).$$ Consider an element
$$r(x_{\overline{\sigma}_m})= (x_{(i_{1,1},i_{2,1})}\cdot \,. \, .\, .\,
\cdot x_{( i_{k_1-1,1},i_{k_1,1})})\cdot \, .\, .\, .\, \cdot
(x_{(i_{1,n},i_{2,n})}\cdot \,. \, .\, .\, \cdot x_{(
i_{k_n-1,n},i_{k_n,n})})\in S_{T_d},$$ where $r$ is the regenerating
homomorphism. By Lemma \ref{hhhh},
$$h_{d,0}=r(x_{\overline{\sigma}_m})\cdot \overline h_m$$ with $\overline
h_m$ such that $(\mathcal S_d)_{\overline h_m}=\mathcal S_d$.

We have
$$\begin{array}{ll}
s= &  x_{\sigma_m'}\cdot h_{d,0}\cdot \overline s'=
x_{\sigma_m'}\cdot r(x_{\overline{\sigma}_m})\cdot \overline
h_m\cdot \overline s'=\\  &
 r(x_{\overline{\sigma}_m}) \cdot
x_{\sigma_{m,0}}\cdot \overline h_m\cdot \overline s'=
x_{\sigma_{m,0}}\cdot r(x_{\overline{\sigma}_m'}) \cdot \overline
h_m\cdot \overline s',
\end{array}$$
where
$x_{\overline{\sigma}_m'}=\lambda(\sigma_{m,0})(x_{\overline{\sigma}_m})$.
We have  $\overline s_1'=r(x_{\overline{\sigma}_m'}) \cdot \overline
h_m\cdot \overline s'\in S_{T_d}$, its length $ln(\overline
s_1')=k$, its image $\alpha(\overline
s_1')=\sigma_{m,0}^{-1}\alpha(s)$, and $(\mathcal S_d)_{\overline
s_1'}=\mathcal S_d$. Therefore, by Theorem \ref{fuc}, $\overline
s_1'=r(x_{\sigma})\cdot h_{d,g}$, where $\sigma=\alpha(\overline
s_1')= \sigma_{m,0}^{-1}\alpha(s)$ and
$g=\frac{k-ln_t(x_{\sigma})}{2}-d+1$.

Now, assume that Theorem \ref{fuc1} is proved for all $m<m_0$ and
consider an element
$$
s=x_{\sigma_1}\cdot \, .\, .\, .\, \cdot x_{\sigma_{m_0}}\cdot
\overline s_1,$$ where $\overline s_1\in S_{T_d}$ has the length
$k\geq 3(d-1)$ and it is such that $(\mathcal S_d)_{\overline
s_1}=\mathcal S_d$. We have
$$\begin{array}{l}
s=x_{\sigma_1}\cdot \, .\, .\, .\, \cdot x_{\sigma_{m_0}}\cdot
\overline s_1=  x_{\sigma_2'}\cdot \, .\, .\, .\, \cdot
x_{\sigma_{m_0}'}\cdot x_{\sigma_1}\cdot \overline s_1=
\\
x_{\sigma_2'}\cdot \, .\, .\, .\, \cdot x_{\sigma_{m_0}'}\cdot
x_{\sigma_{1,0}}\cdot \overline s_1'=  x_{\sigma_{1,0}}\cdot
x_{\sigma_2''}\cdot\, .\, .\, .\, \cdot x_{\sigma_{m_0}''}\cdot
\overline s_1',
\end{array}$$
where $\sigma_j'=\sigma_1\sigma_j\sigma_1^{-1}$ and
$\sigma_j''=\sigma_{1,0}^{-1}\sigma_j'\sigma_{1,0}$ for $j=2,\dots,
m$, and the element $\overline s_1'\in S_d$ is such that
$ln(\overline s_1')=k$ and $(\mathcal S_d)_{\overline s_1'}=\mathcal
S_d$. Therefore, by inductive assumptions, we have
$$s=  x_{\sigma_{1,0}}\cdot
(x_{\sigma_2''}\cdot\, .\, .\, .\, \cdot x_{\sigma_{m_0}''}\cdot
\overline s_1')= x_{\sigma_{1,0}}\cdot (x_{\sigma_{2,0}}\cdot\, .\,
.\, .\, \cdot x_{\sigma_{m_0,0}}\cdot \overline s_1''),$$ where the
element  $\overline s_1''\in S_d$ is such that $ln(\overline
s_1'')=k$ and $(\mathcal S_d)_{\overline s_1''}=\mathcal S_d$. By
Theorem \ref{fuc}, we have $\overline s_1''=r(x_{\sigma})\cdot
h_{d,g}$, where $\sigma =\alpha(\overline s_1'')=(\sigma_{1,0}\dots
\sigma_{m,0})^{-1}\alpha(s)$ and
$g=\frac{k-ln_t(x_{\sigma})}{2}-d+1$.
 \qed

\begin{cor}\label{cor1}  Let $s_i=x_{\sigma_{1,i}}\cdot \, .\, .\,
.\, \cdot x_{\sigma_{m,i}}\cdot \overline s_i$, $i=1,2$, be two
elements of $\Sigma_{d}^{\mathcal S_d}$, where $\overline s_i\in
S_{T_d}$ of length $ln(\overline s_1)=ln(\overline s_2)=k$. Assume
also that $\alpha(s_1)=\alpha(s_2)$ and $\tau(s_1)=\tau(s_2)$. If
$k\geq 3(d-1)$, then $s_1=s_2$.
\end{cor}
\begin{cor} The Hurwitz element $h_{d,[\frac{d}{2}]}$ is a
stabilizing element of $\Sigma_d$, that is, the  semigroup
$\Sigma_d$ is stable. \end{cor}

\subsection{Factorizations in $\mathcal S_3$} Consider the semigroups $\Sigma_{3,\bf{1}}\subset \Sigma_3$. The semigroup
$\Sigma_3$ is generated by the elements $x_{(1,2)}$, $x_{(1,3)}$,
$x_{(2,3)}$, $x_{(1,2,3)}$, and $x_{(1,3,2)}$ satisfying the
following relations:

\begin{equation} \label{eq1}
x_{(1,2)}\cdot x_{(1,3)}=x_{(2,3)}\cdot x_{(1,2)}=x_{(1,3)}\cdot
x_{(2,3)};\end{equation}
\begin{equation} \label{eq2}
x_{(1,3)}\cdot x_{(1,2)}=x_{(2,3)}\cdot x_{(1,3)}=x_{(1,2)}\cdot
x_{(2,3)};\end{equation}
\begin{equation} \label{eq3}  x_{(1,2)}\cdot x_{(1,2,3)}=x_{(1,3,2)}\cdot x_{(1,2)}=x_{(2,3)}\cdot
x_{(1,3,2)}=x_{(1,2,3)}\cdot x_{(2,3)} ;\end{equation}
\begin{equation} \label{eq4}
x_{(1,2)}\cdot x_{(1,3,2)}=x_{(1,2,3)}\cdot x_{(1,2)}=x_{(1,3)}\cdot
x_{(1,2,3)}=x_{(1,3,2)}\cdot x_{(1,3)} ;\end{equation}
\begin{equation} \label{eq5}  x_{(2,3)}\cdot x_{(1,2,3)}=x_{(1,3,2)}\cdot x_{(2,3)}=x_{(1,3)}\cdot
x_{(1,3,2)}=x_{(1,2,3)}\cdot x_{(1,3)} ;\end{equation}
\begin{equation} \label{eq6}
x_{(1,3)}\cdot x_{(1,3,2)}=x_{(1,2,3)}\cdot x_{(1,3)}=x_{(2,3)}\cdot
x_{(1,2,3)}=x_{(1,3,2)}\cdot x_{(2,3)}, \end{equation}
\begin{equation} \label{eq7}
x_{(1,2,3)}\cdot x_{(1,3,2)}=x_{(1,3,2)}\cdot x_{(1,2,3)}.
\end{equation}

Denote by $$\begin{array}{l} s_1=x_{(1,2)}^2, \quad s_2=x_{(2,3)}^2,
\quad s_3=x_{(1,3)}^2, \quad s_4=x_{(1,2,3)}\cdot x_{(1,3,2)}, \\
s_5=x_{(1,2,3)}\cdot x_{(1,2)}\cdot x_{(2,3)}, \quad
s_6=x_{(1,2,3)}^3, \quad s_7=x_{(1,3,2)}^3. \end{array}$$ It is easy
to see that $s_1,\dots, s_7\in \Sigma_{3,\bf{1}}$.

\begin{thm} \label{S3} The semigroup $\Sigma_{3,\bf{1}}$ has the following
presentation:
$$\begin{array}{ll}
\Sigma_{3,\bf{1}}=  \{ s_1,\dots,s_7\mid   & s_i\cdot s_j=s_j\cdot
s_i \quad \text{for}\quad 1\leq i,j\leq 7; \\ &  s_i\cdot s_k=
s_j\cdot s_k \quad \text{for}\,\, 1\leq i,j\leq 3,\,\, 4\leq k\leq
7;
\\ &  s_i\cdot s_6=
s_i\cdot s_7 \quad \text{for}\,\, 1\leq i\leq 3;
\\ & s_1\cdot s_2=s_1\cdot s_3=s_2\cdot s_3;
 \\  &
 s_4^3=s_6\cdot s_7; \\ &  s_5^2=s_1^2\cdot s_4  \quad s_5^3=s_1^3\cdot s_6; \\ &
 s_4\cdot s_5=s_1\cdot s_6=s_1\cdot s_7
 \} .
\end{array} $$
\end{thm}
\proof First of all let us show that the elements $s_1,\dots, s_7$
generate $\Sigma_{3,\bf{1}}$. Indeed, assume that any
$s\in\Sigma_{3,\bf{1}}$ of length $ln(s)\leq k$ can be written as a
word in $s_1,\dots , s_7$ and consider an element
$s\in\Sigma_{3,\bf{1}}$ of length $ln(s)=k+1$. Moving the factors
$x_{(1,2,3)}$ and $x_{(1,3,2)}$ to the left side, any element $s\in
\Sigma_{3,\bf{1}}$ can be written in the following form
$$s=x_{(1,2,3)}^a\cdot x_{(1,3,2)}^b\cdot s',$$ where $a,b$ are non-negative
integers and $s'$ is a word in letters $x_{(1,2)}$, $x_{(1,3)}$, and
$x_{(2,3)}$.

By Lemmas \ref{chain1} and \ref{chain2}, if $ln(s')\geq 3$, then
$s'$ can be written in the form $s'=x_{(i,j)}^2\cdot s''$.
Similarly, if either $a\geq 3$, or $b\geq 3$, or both $a$ and $b$
are positive, then $s=s_i\cdot \widetilde s$, where $i$ is either
$6$, or $7$, or $4$ and $\widetilde s\in \Sigma_{3,\bf{1}}$,
$ln(\widetilde s)\leq k-1$. So we need to consider only the cases
when $ln(s')\leq 2$ and either $0\leq a\leq 2$, $b=0$ or $a=0$,
$0\leq b\leq 2$. If $a=b=0$, then it is obvious that $s'=s_i$ for
some $i=1,2,3$, since $s=s'\in \Sigma_{3,\bf{1}}$.

Consider the case $a=1$ and $b=0$, that is, $s=x_{(1,2,3)}\cdot s'$.
Since $s\in \Sigma_{3,\bf{1}}$ and $\alpha(x_{(1,2,3)})=(1,2,3)$, we
have $\alpha(s')=(1,3,2)$. Therefore $s'$ is equal to either
$x_{(1,2)}\cdot x_{(2,3)}$, or $x_{(1,3)}\cdot x_{(1,2)}$, or
$x_{(2,3)}\cdot x_{(1,3)}$. But, by (\ref{eq2}), the last three
elements are equal to each other and in this case $s=s_5$.

Similarly, if $a=0$, $b=1$, that is, $s=x_{(1,3,2)}\cdot s'$, then
we obtain that  $s'$ is equal to either $x_{(1,3)}\cdot x_{(2,3)}$,
or $x_{(2,3)}\cdot x_{(1,2)}$, or $x_{(1,2)}\cdot x_{(1,3)}$, and,
by (\ref{eq1}), the last three elements are equal to each other.
Therefore, by (\ref{eq4}), we have
$$s=x_{(1,3,2)}\cdot x_{(1,3)}\cdot x_{(2,3)}=x_{(1,3)}\cdot
x_{(1,2,3)}\cdot x_{(2,3)}=x_{(1,2,3)}\cdot x_{(1,2)}\cdot
x_{(2,3)}=s_5.$$

If  $a=2$, $b=0$, that is, $s=x_{(1,2,3)}^2\cdot s'$, then we obtain
that $\alpha (s')=(1,2,3)$ and hence $s'=x_{(2,3)}\cdot x_{(1,2)}$.
Therefore, by (\ref{eq3}), $$s=x_{(1,2,3)}^2\cdot x_{(2,3)}\cdot
x_{(1,2)}=x_{(1,2,3)}\cdot x_{(2,3)}\cdot x_{(1,3,2)}\cdot
x_{(1,2)}=x_{(1,2,3)}\cdot x_{(1,3,2)}\cdot x_{(1,2)}\cdot
x_{(1,2)}=s_4\cdot s_1.$$

Finally, if  $a=0$, $b=2$, that is, $s=x_{(1,3,2)}^2\cdot s'$, then
we have $\alpha (s')=(1,3,2)$ and hence $s'=x_{(1,3)}\cdot
x_{(1,2)}$. Therefore, by (\ref{eq4}),
$$s=x_{(1,3,2)}^2\cdot x_{(1,3)}\cdot x_{(1,2)}=x_{(1,3,2)}\cdot
x_{(1,3)}\cdot x_{(1,2,3)}\cdot x_{(1,2)}=x_{(1,3,2)}\cdot
x_{(1,2,3)}\cdot x_{(1,2)}\cdot x_{(1,2)}=s_4\cdot s_1$$ and as a
result we obtain that $\Sigma_{3,\bf{1}}$ is generated by
$s_1,\dots, s_7$.

Since the inspection, that the generators $s_1,\dots, s_7$ of
$\Sigma_{3,\bf{1}}$ satisfy all relations mentioned in the statement
of Theorem \ref{S3}, is similar, we will check only one of them and
the inspection of all other relations will be left to the reader.

Let us show, for example, that $s_4\cdot s_5=s_6\cdot s_1$. By
(\ref{eq1}) -- (\ref{eq7}), we have
$$\begin{array}{l} s_4\cdot s_5=x_{(1,2,3)}\cdot \underline{x_{(1,3,2)}\cdot x_{(1,2,3)}}\cdot
x_{(1,2)}\cdot x_{(2,3)}=x_{(1,2,3)}\cdot (x_{(1,2,3)}\cdot
\underline{x_{(1,3,2)})\cdot x_{(1,2)}}\cdot x_{(2,3)}= \\
x_{(1,2,3)}\cdot x_{(1,2,3)}\cdot (x_{(1,2)}\cdot
\underline{x_{(1,2,3)})\cdot x_{(2,3)}}=x_{(1,2,3)}\cdot
x_{(1,2,3)}\cdot x_{(1,2)}\cdot (\underline{x_{(1,2)}\cdot x_{(1,2,3)}})= \\
x_{(1,2,3)}\cdot x_{(1,2,3)}\cdot \underline{x_{(1,2)}\cdot
(x_{(1,3,2)}}\cdot x_{(1,2)})=x_{(1,2,3)}\cdot x_{(1,2,3)}\cdot
(x_{(1,2,3)}\cdot x_{(1,2)})\cdot x_{(1,2)}=s_6\cdot s_1.
\end{array}$$

The statement that the relations, mentioned in Theorem \ref{S3}, are
defining follows from the next theorem.

\qed

\begin{thm} \label{S31}
Each element $s\in \Sigma_{3,\bf{1}}$, $s\neq \bf{1}$, has a normal
form, that is, it is equal to one and the only one element of the
following form
$$s=\left\{
\begin{array}{ll} s_i^n, & \quad i=1,2,3, \quad n\in\mathbb N, \\
s_4^a\cdot s_6^m\cdot s_7^n,  & \quad 0\leq a\leq 2,\, m\geq
0,n\geq 0,\, a+m+n>0, \\
s_1^n\cdot s_2, & \quad n\in\mathbb N, \\
s_1^n\cdot s_6^m, &  \quad m,n \in \mathbb N,
\\
s_1^n\cdot s_5\cdot s_6^m, &  \quad m\geq 0, \quad n\geq 0, \\
s_1^n\cdot s_4\cdot s_6^m, &   \quad m\geq 0, \quad n\geq 0.
\end{array} \right. $$
\end{thm}

\proof If $s\not\in \Sigma_{3,\bf{1}}^{\mathcal S_3}$, then it is
obvious that $s$ is equal either $s_i^n$,  $i=1,2,3$, or $s_4^a\cdot
s_6^m\cdot s_7^n$.

Let $s\in \Sigma_{3,\bf{1}}^{\mathcal S_3}$. If $s\in
S_{T_3,\bf{1}}$, then by Clebsch -- Hurwitz Theorem $s=h_{3,g}$ for
some $g$.

Let $s=s'\cdot s''$, where $s'=x_{(1,2,3)}^{k_1}\cdot
x_{(1,3,2)}^{k_2}$ and $s''\in S_{T_3}$. Applying relations
(\ref{eq3}) -- (\ref{eq6}), we can assume that $s'=x_{(1,2,3)}^k$
for $k=k_1+k_2$. If $k\equiv 0\, (mod\, 3)$, then by relations in
Theorem \ref{S3}, we have $s=s_1^n\cdot s_6^m$. If $k\equiv 1\,
(mod\, 3)$, then $s'=s_6^m\cdot x_{(1,2,3)}$ and $x_{(1,2,3)}\cdot
s''\in \Sigma_{3,\bf{1}}$. By Theorem \ref{S3}, $x_{(1,2,3)}\cdot
s''=s_5\cdot s_1^n$ for some $n\geq 0$. Similarly, if $k\equiv 2\,
(mod\, 3)$, then $s'=s_6^m\cdot x_{(1,2,3)}^2$ and
$x_{(1,2,3)}^2\cdot s''\in \Sigma_{3,\bf{1}}$. Applying relations
(\ref{eq3}) -- (\ref{eq6}), we get $x_{(1,2,3)}^2\cdot
s''=x_{(1,2,3)}\cdot x_{(1,3,2)}\cdot s_1''=s_4\cdot s_1''$ for some
$s_1''\in S_{T_3,\bf{1}}$, and by relations in Theorem \ref{S3}, we
obtain that $s=s_1^n\cdot s_4\cdot s_6^m$. \qed

\begin{thm} \label{ttt} Up to simultaneous
conjugation, an element $\overline s \in \Sigma_3$  is equal either
to $s$, where $s$ is an element of $\Sigma_{3,\bf{1}}$ described in
Theorem {\rm \ref{S31}}, or to
$$\overline s=\left\{
\begin{array}{ll} x_{(1,2)}^{2k+1}, &  \, \,
k\geq 0, \\
x_{(1,2,3)}^n\cdot x_{(1,3,2)}^m,  & \, \, n> m,\, \, n\, \,
\text{or}\, \, m\not\equiv 0\, (mod\, 3), \\
x_{(1,2)}^n\cdot x_{(2,3)}, & \, \,
n\in\mathbb N, \\
x_{(1,2)}^{n}\cdot x_{(1,2,3)}^{3m}\cdot x_{(1,3,2)}^a, &  \, \, n
\in \mathbb N,\, m\geq 0,  \, a=0,1,2,  \, \text{and}\, \, a\neq 0
\, \text{if}\, \, n\equiv 0(mod\, 2).
\end{array} \right. $$
\end{thm}

\proof To prove Theorem \ref{ttt}, one must consider separately the
following cases:

1)$(\mathcal S_3)_s=\mathcal S_2$;

2) $(\mathcal S_3)_s=A_3$, where $A_3$ is the alternating group;

3) $s\in S_{T_3}$, $(\mathcal S_3)_s=\mathcal S_3$, and $\alpha(s)$
is either a transposition or a cyclic permutation of length $3$;

4) $s\not\in S_{T_3}$, $(\mathcal S_3)_s=\mathcal S_3$,  and
$\alpha(s)$ is either a transposition or a cyclic permutation of
length $3$.
\newline It is easy to see that in the first three casees $s$ is equal (up to conjugation) respectively to
1) $x_{(1,2)}^{2k+1}$; 2) $x_{(1,2,3)}^n\cdot x_{(1,3,2)}^m$, 3)
$x_{(1,2)}^n\cdot x_{(2,3)}$.

In case 4) we have $s=s_1\cdot s_2$, $s_1\in S_{T_d}$ and $s_2$ is
represented as a word in letters $x_{(1,2,3)}$ and $x_{(1,3,2)}$. By
(\ref{eq3}) and (\ref{eq4}), we can assume that $s_1= x_{(1,2)}^n$.
Next, we have $$x_{(1,2)}\cdot x_{(1,2,3)}^3= x_{(1,3,2)}^3\cdot
x_{(1,2)}=x_{(1,2)}\cdot x_{(1,3,2)}^3.$$ Applying these relations
and (\ref{eq7}), we obtain that $s=x_{(1,2)}^{n}\cdot
x_{\sigma}^{3m}\cdot x_{\sigma^{-1}}^a$, where $\sigma=(1,2,3)$ or
$(1,3,2)$. To complete the proof, notice that
$\lambda((1,2))(x_{\sigma})=x_{\sigma^{-1}}$. \qed

\begin{cor} \label{uniq} Let $(\mathcal S_3)_s=\mathcal S_2$ or $\mathcal S_3$ for $s \in \Sigma_3$.
Then $s$ is uniquely defined up to simultaneous conjugation by its
type $\tau(s)$ and the type $t(\alpha(s))$ of its image
$\alpha(s)\in \mathcal S_3$.

Up to simultaneous conjugation, there are exactly $[\frac{n}{6}]+1$
different elements $s\in \Sigma_{3,\bf{1}}^{A_3}$ of $ln(s)=n$, and
if $\alpha(s)\neq\bf{1}$, then there are exactly $m=-[\frac{-n}{3}]$
different elements $s\in \Sigma_3^{A_3}$ of $ln(s)=n$.
\end{cor}
\subsection{Cayley's imbeddings}
As is well-known, any finite group $G$ can be embedded into some
symmetric group. In particular, if $N=|G|$ is the order of a group
$G$, then we can have Cayley's imbedding $c:G\hookrightarrow
Sym(G)\simeq \mathcal S_N$:
$$(g_1)\sigma_g=g_1g\quad \text{for}\, g,g_1\in G,\,\,
c(g)=\sigma_g,$$ that is, $G$ acts on itself by multiplication from
the right side. Let us identify the group $G$ with its image $c(G)$
and denote by $N(G)$ and $C(G)$ the normalizer and centralizer of
$G$ in $\mathcal S_N$, respectively. Since $N(G)$ acts on $G$ by
conjugations, we have the natural homomorphism $a:N(G)\to Aut(G)$.
\begin{thm} \label{aut} Let $c:G\hookrightarrow Sym(G)\simeq \mathcal S_N$ be
the Cayley's imbedding of a finite group $G$. Then the natural
homomorphism $a:N(G)\to Aut(G)$ has the following properties:
\begin{itemize}
\item[($i$)] $a$ is an epimorphism,
\item[($ii$)] $\ker a=C(G)\simeq G$,
\item[($iii$)] the group generated by $G$ and $C(G)$ is isomorphic
to the amalgamated direct product $G\times_C G$, where $C$ is the
center of $G$.
\end{itemize}
\end{thm}
\proof Consider an automorphism $f\in Aut(G)$ as a permutation
$\sigma_f\in \mathcal S_N$ of the elements of $G$:
$$(g)\sigma_f=f(g)\, \, \text{for}\, \, g\in G.$$
Let us show that $\sigma_f\in N(G)$. For all $g_1\in G$ we have
$$ (g_1)\sigma_f^{-1} \sigma_g\sigma_f=(f^{-1}(g_1))\sigma_g\sigma_f=(f^{-1}(g_1)g)\sigma_f=
f(f^{-1}(g_1)g)=g_1f(g)=(g_1)\sigma_{f(g)},$$ that is,
$\sigma_f^{-1} \sigma_g\sigma_f=\sigma_{f(g)}\in G$ for all $g\in
G$. Hence $\sigma_f\in N(G)$ and, moreover, the conjugation of the
elements of $G$ by $\sigma_f$ defines the automorphism $f$ of the
group $G$. Therefore the homomorphism $a$ is an epimorphism.

It is obvious that $C(G)=\ker a$. Consider $\sigma\in C(G)$. We have
$\sigma_g \sigma=\sigma \sigma_g$ for all $g\in G$. Therefore
$$(g_1)\sigma_g \sigma=(g_1g)\sigma=((g_1)\sigma)\cdot g$$
for all $g_1, g\in G$. In particular, for $g_1=\bf{1}$ if we denote
$({\bf{1}})\sigma$ by $g_{\sigma}$, then we have
$$({\bf{1}} )\sigma_g \sigma=(g)\sigma=g_{\sigma} g$$
for all $g\in G$. The equality $(g)\sigma=g_{\sigma}g$ shows that
$\sigma$ acts on $G$ as multiplication in $G$ from the left side by
the element $g_{\sigma}\in G$. Obviously, the multiplications by
elements of $G$ from the left side and from the right side commute.
Therefore $C(G)\simeq G$.

Remind that, by definition, the group $G$ acts on itself by the
multiplication from the right side. It is easy to see from this that
the group generated by $G$ and $C(G)$ is isomorphic to the
amalgamated direct product $G\times_C G$, where $C$ is the center of
$G$. \qed \\

Any imbedding $G\hookrightarrow \mathcal S_d$ defines an imbedding
$S(G,O)\hookrightarrow \Sigma_d$. Let $c:S_G=S(G,G)\hookrightarrow
\Sigma_d$ be the imbedding of semigroups defined by Cayley's
imbedding $c:G\to \mathcal S_N$. Theorem \ref{aut} implies the
following
\begin{cor} \label{norbits} The orbits of conjugation
action of $\mathcal S_N$ on $\Sigma_N$ intersecting $S(G,G)$ are in
one to one correspondence with the orbits of the action $Aut(G)$ on
$S(G,G)$.
\end{cor}

\section{Hurwitz spaces}

\subsection{Marked Riemannian surfaces}
Let $f:C\to D_R=\{ z\in \mathbb C\mid |z|\leq R\}$ be a Riemannian
surface, that is, $f$ is a finite proper continuous ramified
covering of the disc $D_R=\{ |z|\leq R\}$ (or $\mathbb P^1$ if
$R=\infty$) of degree $d$ branched at finite number of points in
$D_R^0=D_R\setminus
\partial D_R=\{ |z|<R\}$ (it is not assumed that $C$ is necessary to be connected).
Two coverings $(C^{\prime},f^{\prime})$ and
$(C^{\prime\prime},f^{\prime\prime})$ of $D_R$ are said to be
isomorphic if there is a homeomorphism $h:C^{\prime}\to
C^{\prime\prime}$ preserving the orientation and such that
$f^{\prime}=h\circ f^{\prime\prime}$, and they are said to be {\it
equivalent} if there are preserving orientations homeomorphisms
$\psi :D_R\to D_R$ and $\varphi : C^{\prime}\to C^{\prime\prime}$
such that $\psi$ leaves fixed the boundary $\partial D_R$ and
$\psi\circ f^{\prime}=f^{\prime\prime}\circ\varphi$. Denote by
${\mathcal{R}}_{R,d}$ the set of equivalence classes of the
coverings of $D_R$ of degree $d$ with respect to this equivalence.

Let $q_1, \dots, q_b\in D_R^0$ be the points over which $f$ is
ramified. Let us fix the point $o=o_R=e^{\frac{3}{2}\pi i}R\in
\partial D_R$ (if $R=\infty$, then, by definition, $o_{\infty}=\infty=\mathbb P^1\setminus \mathbb C$)
and number the points in $f^{-1}(o)$. A numbering of the points in
$f^{-1}(o)$ defines an order on the points in  $f^{-1}(o)$. Such
coverings $(C,f)$ with fixed point $o\in D_R$ and fixed ordering of
the points of $f^{-1}(o)$ will be called coverings with {\it ordered
set of sheets} or a {\it marked coverings}. We say that marked
coverings $(C^{\prime},f^{\prime})_m$ and
$(C^{\prime\prime},f^{\prime\prime})_m$ are {\it equivalent} if
there are homeomorphisms  $\psi :D_R\to D_R$ and $\varphi :
C^{\prime}\to C^{\prime\prime}$ preserving  orientations and such
that
\begin{itemize}
\item[($i$)] $\psi$ leaves fixed the boundary $\partial
D_R$;  \item[($ii$)] $\varphi(p_i^{\prime})=p_i^{\prime\prime}\in
f^{\prime\prime^{-1}}(o)$ for each $p_i^{\prime}\in
f^{\prime^{-1}}(o)$, $i=1,\dots,d$;
\item[($iii$)] $\psi\circ f^{\prime}=f^{\prime\prime}\circ\varphi$.
\end{itemize}
Denote by ${\mathcal{R}}_{R,d}^m$ the set of equivalence classes of
the marked coverings of $D_R$ of degree $d$ with respect to this
equivalence. Renumberings of sheets define an action of the
symmetric group $\mathcal S_d$ on ${\mathcal{R}}_{R,d}^m$ and it is
easy to see that ${\mathcal{R}}_{R,d}={\mathcal{R}}_{R,d}^m/\mathcal
S_d$.

If $R_1<R_2<\infty$, then any ramified covering $f:C\to D_{R_1}$ can
be extended to a ramified covering $\tilde f:\tilde C\to D_{R_2}$
non-ramified over $D_{R_2}\setminus D_{R_1}$. The lift of the path
$$l(t)=e^{\frac{3}{2}\pi i}(R_2t+(1-t)R_1)\subset D_{R_2}\setminus
D^0_{R_1},\qquad t\in [0,1],$$ to $\tilde C$ defines $d$ paths
$\tilde f^{-1}(l(t))$ connecting the points of $f^{-1}(o_{R_1})$ and
$f^{-1}(o_{R_2})$. If $(C,f)_m$ is a marked covering, then these
paths transfer the order from the set $f^{-1}(o_{R_1})$ to the set
$f^{-1}(o_{R_2})$. As a result, we obtain an isomorphism
$i_{R_1,R_2}:{\mathcal{R}}_{R_1,d}^m \hookrightarrow
{\mathcal{R}}_{R_2,d}^m$.

Similarly, for any marked covering $(C,f)_m$ of $\mathbb P^1$ and
for any $R>0$ there is an equivalent covering $(\overline
C,\overline f)_m$ those branch points belong to $D_R^0$. Consider
the restriction $\tilde f$ of $\overline f$ to $\tilde C=\overline
f^{-1}(D_R)$. If we lift the path
$$l(t)=e^{\frac{3}{2}\pi i}R/t\subset \mathbb P^1\setminus D_R^0,\qquad t\in [0,1],$$ to
$\overline C$, then we obtain $d$ paths $\overline f^{-1}(l(t))$
connecting the points of $f^{-1}(o_{\infty})$ and $f^{-1}(o_R)$
which transfer the order from $\overline f^{-1}(o_{\infty})$ to the
set $\tilde f^{-1}(o_R)$. Obviously, the equivalence class of
obtained marked covering $(\tilde C,\tilde f)_m$ does not depends on
the choice of a representative $(\overline C,\overline f)_m$.
Therefore we obtain an imbedding of $i_{\infty
,R}:{\mathcal{R}}_{\infty ,d}^m \hookrightarrow
{\mathcal{R}}_{R,d}^m$. It is easy to see that $i_{\infty
,R_2}=i_{R_1,R_2}\circ i_{\infty ,R_1}$ for any $R_2\geq R_1>0$.

\subsection{Semigroups of marked coverings} A closed loop $\gamma\subset D_R\smallsetminus
\{q_1,\dots,q_b\}$ starting and ending at $o=o_R$ can be lifted to
$C$ by means of $f$ and we get $d$ paths staring and ending at the
points in $f^{-1}(o)$. Such lift of the loops defines a homomorphism
(the {\it monodromy of marked covering}) $\mu
:\pi_1(D_R\smallsetminus \{q_1,\dots,q_b\}, o)\to \mathcal S_d$ to
the symmetric group $\mathcal S_d$ (the monodromy sends starting
points of the lifted paths to the ends of the corresponding paths).
Conversely, if a homomorphism $\mu :\pi_1(D_R\smallsetminus
\{q_1,\dots,q_b\}, o)\to \mathcal S_d$ is given, then it defines a
marked covering  $f:C\to D$ whose monodromy is $\mu$.

The fundamental group $\pi_1(D_R\smallsetminus \{ q_1,\dots,q_b\},
o)$ is generated by loops $\gamma _1,\dots,\gamma_b$ of the
following form. Each loop $\gamma _i$ consists of a path $l_i$
starting at $o$ and ending at a point $q'_i$ close to $q_i$,
followed by a circuit in positive direction (with respect to the
complex orientation on $\mathbb C$) around a circle $\Gamma _i$ of
small radius with the center at $q_i$, $q'_i\in \Gamma$, followed by
the return to $q_0$ along the path $l_i$ in the opposite direction;
for $i\neq j$ the loops $\gamma_i$ and $\gamma_j$ have the only one
common point, namely, $o$;  and the product $\gamma_1\dots \gamma_b
=\partial D_R$ in $\pi_1(D_R\smallsetminus \{ q_1,\dots,q_b\}, o)$.
Such collection of generators is called a {\it good geometric base}
of the group $\pi_1(D_R\smallsetminus \{ q_1,\dots,q_b\}, o)$. It is
well known that if $R<\infty$, then $\gamma _1,\dots,\gamma_b$ are
free generators of $\pi_1(D_R\smallsetminus \{ q_1,\dots,q_b\}, o)$,
that is, $\pi_1(D_R\smallsetminus \{ q_1,\dots,q_b\}, o)= \langle
\gamma _1,\dots,\gamma_b\rangle$; and if $R=\infty$, then $\gamma
_1,\dots,\gamma_b$ generate the group $\pi_1(\mathbb
P^1\smallsetminus \{ q_1,\dots,q_b\}, o)$ being subject to the
relation $\gamma_1\dots\gamma_b=\bf{1}$.

If we choose a good geometric base  $\gamma_1,\dots ,\gamma_b$, then
the monodromy $\mu$ is defined by a collection of elements
$\sigma_1=\mu(\gamma_1),\dots ,\sigma_n=\mu(\gamma_b)\in \mathcal
S_d$ called {\it local monodromies} and the product
$\sigma=\sigma_1\dots \sigma_b=\mu(\partial D)$ is called the {\it
global monodromy} of $f$. It is easy to see that if $R=\infty$, then
the global monodromy is equal to $\bf{1}$.

The collection $( \sigma_1,\dots ,\sigma_b)$  depends on the choice
of a good geometric base $\gamma_1,\dots ,\gamma_b$. Any good
geometric base can be obtained from $\gamma_1,\dots ,\gamma_b$  by
means of a finite sequence of Hurwitz moves. In the other words, the
braid group $\text{Br}_b$ naturally acts on the set of good
geometric bases of $\pi_1(D_R\smallsetminus \{ q_1,\dots,q_b\}, o)$
as the Hurwitz moves (\cite{M-T}). Therefore if $(\sigma'_1,\dots
,\sigma'_b)$ is a collection corresponding to some other good
geometric base $\gamma'_1,\dots,\gamma'_b$, then the collection
$(\sigma'_1,\dots ,\sigma'_b)$ can be obtained from
$(\sigma_1,\dots,\sigma_b)$ by means of a finite sequence of Hurwitz
moves (see subsection \ref{mu}).

Let $R<\infty$. One can define a structure of semigroup on the set
$\mathcal R_{R,d}^m$ as follows. Let $(C_{1},f_1)_m$ and
$(C_{2},f_2)_m$ be two marked coverings of degree $d$. Let us choose
two continuous preserving the orientations imbeddings
$\varphi_j:D_R\to D_R$, $j=1,2$, of the disc $D_R$ to itself leaving
fixed the point $o$ and such that
\begin{itemize}  \item[($i$)] the image
$\varphi_1(D_R)=\{ u\in D_R\mid \text{Re}\, u\geq 0\}$ is the right
halfdisc and
\\ $\varphi_1(\{ u\in
\partial D_R\mid \text{Re}\, u\leq 0\})= \{ u\in D_R\mid \text{Re}\, u=0\}$ is the vertical diameter;
\item[($ii$)] $\varphi_2(D_R)=\{ u\in D_R\mid \text{Re}\, u\leq 0\}$
is the left halfdisc and \\ $\varphi_2(\{ u\in
\partial D_R\mid \text{Re}\, u\geq 0\})= \{ u\in D_R\mid \text{Re}\, u=0\}$.
\end{itemize}
Let us identify the points belonging to the sets $f_1^{-1}(o)$ and
$f_2^{-1}(o)$ by means of the orders on the sets of these points,
and after that let us identify, by continuity, the points belonging
to the $d$ paths $f_1^{-1}(\{ u\in
\partial D_R\mid \text{Re}\, u\leq 0\})$  in $C_{1}$ with the points belonging to the $d$ paths
$f_2^{-1}(\{ u\in
\partial D_R\mid \text{Re}\, u\geq 0\})$ in $C_{2}$ so that the images under the mappings  $\varphi_1\circ f_1$ Х
$\varphi_2\circ f_2$ of the all identified points should be
coincided. By means of this identification, we can glue the surfaces
$C_{1}$ and $C_{2}$ along these $d$ paths and, as a result we obtain
a marked covering $(C,f)_m$, where $f(q)=\varphi_1(f_1(q))$ if $q\in
C_{1}$ and $f(q)=\varphi_2( f_2(q))$ if $q\in C_{2}$. We call the
obtained covering $(C,f)_m$ the {\it product} of marked coverings
$(C_{1},f_1)_m$ and $(C_{2},f_2)_m$ (notation:
$(C,f)_m=(C_{1},f_1)_m\cdot(C_{2},f_2)_m$). It is easy to see that
the product introduced above defines a structure of non-commutative
semigroup on $\mathcal R_{R,d}^m$ such that the maps $i_{R_1,R_2}$
are isomorphisms of semigroups for all $R_1\geq R_2>0$.

It is obvious that the semigroup $\mathcal R_d^m= \mathcal
R_{R,d}^m$ is generated by the marked coverings $(C,f)_m$ which are
coverings of the disc $D=D_R$ with a single branch point $q_1$. Such
coverings are defined uniquely (up to equivalence) by their global
monodromy $\sigma_f=\mu(\partial D)\in \mathcal S_d$ where
$\mu=\mu_f$ is the monodromy of the marked covering $(C,f)_m$.
Therefore the number of generators is equal to $d!$. Denote by
$x_{\sigma_f}$ the generator of the semigroup $\mathcal R_d$
corresponding to a covering $(C,f)_m$ with single branch point. A
simple inspection shows that in the semigroup $\mathcal R_d^m$ the
generators $x_{\sigma}$ satisfy the following defining relations:
$$x_{\sigma_1}\cdot
x_{\sigma_2}=x_{\sigma_2}\cdot\,x_{(\sigma_2^{-1}\sigma_1\sigma_2)},\qquad
x_{\sigma_1}\cdot
x_{\sigma_2}=x_{(\sigma_1\sigma_2\sigma_1^{-1})}\cdot
x_{\sigma_1},$$  and  $x_{\sigma_1}\cdot x_{\bold 1}=x_{\sigma_1}$,
$x_{\bold 1}\cdot x_{\sigma_2}=x_{\sigma_2}$ for all
$\sigma_1,\sigma_2\in \mathcal S_d$.

It is easy to check that if a marked covering $(C,f)_m$ is equal to
$x_{\sigma_1}\cdot .\, .\, .\cdot x_{\sigma_n}$ in $\mathcal R_d^m$,
then its global monodromy $\sigma_f=\mu(\partial D)$ is equal to the
product $\sigma_1\dots\sigma_n$ and it is obvious that the
comparison to each marked  covering its global monodromy defines a
homomorphism from $\mathcal R_d^m$ to the symmetric group $\mathcal
S_d$. Denote this homomorphism by $\alpha :\mathcal R_d^m\to
\mathcal S_d$.

Renumberings of the sheets of the marked coverings define an action
of $\mathcal S_d$ on $\mathcal R_d^m$. Namely, an element
$\sigma_0\in\mathcal S_d$ acts on the generators $x_{\sigma}$ by the
following rule: $x_{\sigma}\mapsto
x_{(\sigma_0^{-1}\sigma\sigma_0)}$. This action defines a
homomorphism $\lambda:\mathcal S_d\to \text{Aut}(\mathcal R_d^m)$.
Therefore we obtain the following
\begin{prop} \label{ppp} The semigroup $\mathcal R_d^m$ as a semigroup over $\mathcal S_d$
is naturally isomorphic to $\Sigma_d$. \end{prop}

According to Proposition \ref{ppp}, we call the elements of
$\Sigma_d$ {\it monodromy factorizations} of the coverings of degree
$d$.

It is easy to see that the kernel $\ker \alpha=\mathcal R_{d,\bold
1}^m=\{ (C,f)_m\in \mathcal R_d^m\mid \sigma_f=\bold 1\}$ is a
subsemigroup in $\mathcal R_d^m$ isomorphic to $\Sigma_{d,\bf{1}}$
and if the disc $D$ is embedded in $\mathbb P^1$, then the elements
of $\mathcal R_{d,\bold 1}^m$ are the marked coverings $f:C\to D$
for which there are extensions to marked coverings $\widetilde
f:\widetilde C\to \mathbb C\mathbb P^1$ non-ramified over $\mathbb
P^1\setminus D$. Note that the extension $\widetilde f:\widetilde
C\to \mathbb C\mathbb P^1$ of a marked covering $f:C\to D$ with the
global monodromy $\mu_f(\partial D)=\bf{1}$ is defined uniquely up
to equivalence.

The inverse statement is also true: the image of $\mathcal
R_{\infty,d}^m$ under the imbedding $i_{\infty, R}$ coincides with
$\mathcal R_{d,\bold 1}^m$. In the sequel we will identify $\mathcal
R_{\infty,d}^m$ with the semigroup $\mathcal R_{d,\bold 1}^m$ by
means of this isomorphism. As a result, we have the following
\begin{prop} \label{res} On the set of equivalence classes of marked coverings of
$\mathbb P^1$ of degree $d$ there is a natural semigroup structure
isomorphic to $\Sigma_{d,{\bf{1}}}$. \end{prop}

\subsection{Hurwitz spaces of marked Riemannian surfaces}
In this subsection we describe the Hurwitz spaces
$\text{HUR}_d^m(D)$ of marked ramified degree $d$ coverings of
$D=D_R$ considered up to isomorphisms. The space
$\text{HUR}_d^m(D)=\bigsqcup_{b=0}^{\infty}\text{HUR}_{d,b}^m(D)$ is
the disjoint union of the spaces of coverings branched at $b$
points, $b\in \mathbb N$.

As in \cite{F}, let us consider the symmetric product $D^{(b)}$ of
$b$ copied of $D^0=D\setminus \partial D$. It is a complex manifold
of dimension $b$ obtained as the quotient of the cartesian product
$D^b=D^0\times\dots\times D^0$ (with $b$ factors) under the action
of $\mathcal S_b$ which permutes the factors. The points of
$D^{(b)}$ will be identified with the sets of unordered $b$-tuples
of points of $D^0$. Those $b$-tuples which contain fewer than $b$
distinct points form the {\it discriminant locus} $\Delta$ of
$D^{(b)}$.

For a point $B_0=\{ q_{1,0},\dots, q_{b,0}\}\in D^{(b)}\setminus
\Delta$ let us fix the ordered subset $B_0=\{ q_{1,0},\dots,
q_{b,0}\}\subset D$ and choose a good geometric base
$\gamma_1,\dots, \gamma_b$ of $\pi_1(D\setminus B_0,o)$. Then any
word $w$ of the set of words $W_b$ of length $b$ in the letters
$x_{\sigma}$, $\sigma\in \mathcal S_d$, defines a marked covering
$f=f_w: C\to D$ branched over $B_0$ and whose monodromy is $\mu$
such that $\mu(\gamma_i)=\sigma_i$, where $x_{\sigma_i}$ is a letter
in $w$ standing at the $i$-th place.

The choice of a good geometric base allow us to choose the standard
generators $a_1, \dots, a_{b-1}$ in $\pi_1(D^{(b)}\setminus \Delta,
B_0)\simeq \text{Br}_b$ so that this choice defines an action of
$\text{Br}_b$ on the set of words $W_b$ (see subsection \ref{mu}).
In the other words, this choice defines a homomorphism
$\theta_{d,b,R}:\pi_1(D^{(b)}\setminus \Delta, B_0)\simeq
\text{Br}_b\to \mathcal S_N$, where $N=(d!)^b$.

The homomorphism $\theta_{d,b,R}$ allows us to define the space
$\text{HUR}_{d,b}^m(D)$ as an unramified covering
$h_{d,b,R}:\text{HUR}_{d,b}^m(D)\to D^{(b)}\setminus \Delta$
associated with $\theta_{d,b,R}$. Indeed, if we fix a marked
covering $f:C\to D$ with monodromy $\mu$ such that
$\mu(\gamma_i)=\sigma_i$, then any path $\delta(t)$, $0\leq t\leq
1$, in $D^{(b)}$ starting at $B_0$ can be lifted to $D$ and we
obtain $b$ paths $\delta_i(t)$ in $D$ starting at the points
$q_{1,0},\dots, q_{b,0}$. These paths define (up to isotopy) a
continuous family of homeomorphisms $\overline{\delta}_t:D\setminus
B_0\to D\setminus\{\delta_1(t),\dots,\delta_b(t)\}$ leaving fixed
the boundary $\partial D$ such that $\overline{\delta}_0=Id$ and we
can consider a continuous family of marked coverings $f_t:C_t\to D$
branched at $\delta_1(t),\dots,\delta_b(t)$ and given by monodromy
$\mu_t$ such that
$\mu_t(\overline{\delta}_{t*}(\gamma_i))=\sigma_i$. It is obvious
that if $\delta(t)$ is a loop, then the collection
$(\mu_1(\gamma_1),\dots, \mu_1(\gamma_b))$ is Hurwitz equivalent to
$(\mu_0(\gamma_1),\dots, \mu_0(\gamma_b))$. Therefore {\it the
points of the covering space} $\text{HUR}_{d,b}^m(D)$ {\it of the
covering} $h_{d,b,R}:\text{HUR}_{d,b}^m(D)\to D^{(b)}\setminus
\Delta$ {\it naturally parametrize all the marked coverings of $D$
of degree $d$ branched at $b$ points}. The degree of the covering
$h_{d,b,R}$ is equal to $(d!)^b$. As a result, we obtain the
following
\begin{prop} \label{irred1} The irreducible components of $\text{HUR}_{d,b}^m(D)$
are in one to one correspondence with the elements $s$ of the
semigroup $\Sigma_d$ of length $ln(s)=b$.

There is a natural structure of a semigroup on the set of
irreducible components of $\text{HUR}_{d}^m(D)$ isomorphic to
$\mathcal R_d\simeq \Sigma_d$.
\end{prop}

For $R_2\geq R_1>0$ we have the imbedding
$D_{R_1}^{(b)}\hookrightarrow D_{R_2}^{(b)}$ and it is easy to see
that the restriction of $h_{d,b,R_2}$ to
$h_{d,b,R_2}^{-1}(D_{R_1}^{(b)}\setminus \Delta)$ can be identified
with $h_{d,b,R_1}:\text{HUR}_{d,b}^m(D_{R_1})\to
D_{R_1}^{(b)}\setminus \Delta$ by means of $i_{R_1,R_2}$.

According to Proposition \ref{irred1}, we will denote by
$\text{HUR}_{d,s}^m(D)$ the irreducible component of
$\text{HUR}_{d,ln(s)}^m(D)$ corresponding to an element $s\in
\Sigma_d$. In particular, the global monodromy
$\sigma_f=\mu(\partial D)=\alpha(s)\in \mathcal S_d$ is an invariant
of the irreducible component $\text{HUR}_{d,s}^m(D)$. Put
$$\text{HUR}_{d,b,\sigma}^m(D)=\bigcup_{\begin{array}{c} \alpha(s)=\sigma \\ ln(s)=b\end{array}}\text{HUR}_{d,s}^m(D).$$
It follows from consideration above that
$$\text{HUR}_{d,b}^m(\mathbb
P^1)=\bigcup_{R>0}\text{HUR}_{d,b,\bf{1}}^m(D_R).$$

For a fixed type $t$ of elements $s\in \Sigma_d$ let us denote also
by
$$\text{HUR}_{d,t}^m(D)=\bigcup_{\tau(s)=t}\text{HUR}_{d,s}^m(D)$$
and put
$$\text{HUR}_{d,t,\sigma}^m(D)=\text{HUR}_{d,t}^m(D)\cap \text{HUR}_{d,\sigma}^m(D).$$

As it was mentioned above, a marked covering $f:C\to D$ of degree
$d$ branched at the points $q_1,\dots,q_b$ defines (and is defined)
by monodromy $\mu: \pi_1(D\setminus \{ q_1,\dots,q_b\})\to \mathcal
S_d$. The image $\mu(\pi_1(D\setminus \{
q_1,\dots,q_b\}))=Gal(f)\subset \mathcal S_d$ is called the {\it
Galois group} of the covering $f$. It is easy to see that $Gal(f)=
(\mathcal S_d)_s$ if the covering $f$ belongs to
$\text{HUR}_{d,s}^m(D)$. It is not hard to show that the covering
space $C$ of a marked covering $(C,f)_m$ is connected if and only if
the Galois group $Gal(f)$ acts transitively on the set $I_d=[1,d]$.

Denote by $\text{HUR}_{d}^{m,G}(D)$ the union of irreducible
components of $\text{HUR}_{d}^m(D)$ consisting of the coverings with
the Galois group $Gal(f)=G\subset \mathcal S_d$ and put
$\text{HUR}_{d,t}^{m,G}(D)=\text{HUR}_{d}^{m,G}(D)\cap
\text{HUR}_{d,t}^m(D)$ and
$\text{HUR}_{d,t,\sigma}^{m,G}(D)=\text{HUR}_{d,t}^{m,G}(D)\cap
\text{HUR}_{d,t,\sigma}^m(D)$.

By Corollary \ref{cor1}, we have
\begin{thm} Let the type $t$ of monodromy factorization
contains $k$ transpositions. If $k\geq 3(d-1)$ then each irreducible
component of $\text{HUR}_{d,t}^{m,\mathcal S_d}(D)$ is uniquely
defined by the global monodromy $\sigma_f=\mu(\partial D)\in
\mathcal S_d$ of $(C,f)_m$ belonging to this irreducible component.
\end{thm}

\subsection{Hurwitz spaces of (non-marked) coverings of the disc}
\label{disc} To obtain Hurwitz space $\text{HUR}_{d,b}(D)$ of degree
$d$ coverings of a disc $D=D_R$ branched over $b$ points lying in
$D^0$, we must identify all marked coverings of $D$ differ only in
numberings of sheets. The renumberings of sheets induces the action
of $\mathcal S_d$ on the marked fibres. Remind that the actions of
$\text{Br}_b$ and $\mathcal S_d$ on $W_b$ commute. Therefore this
action of $\mathcal S_d$ induces an action on
$\text{HUR}_{d,b}^m(D)$ and we obtain that the space
$\text{HUR}_{d,b}(D)$ is the quotient space:
$\text{HUR}_{d,b}(D)=\text{HUR}_{d,b}^m(D)/\mathcal S_d$. From this
it follows
\begin{prop} \label{irred2} The irreducible components of {\rm
$\text{HUR}_{d,b}(D)$} are in one to one correspondence with the
orbits of the action of $\mathcal S_d$ by simultaneous conjugation
on $\Sigma_{d,b}=\{ s\in \Sigma_d\mid ln(s)=b\}$.
\end{prop}

If $f:C\to D$ is a non-marked covering, then we can also define the
Galois group as $Gal(f)= (\mathcal S_d)_s$. But in this case the
subgroup $Gal(f)\subset\mathcal S_d$ is defined uniquely only up to
inner automorphisms of $\mathcal S_d$.

In the sequel we denote by $\text{HUR}_{\cdot,\cdot,\cdot}(D)$
(resp., $\text{HUR}_{\cdot,\cdot,\cdot}^G(D)$) the image of
introduced above subspaces $\text{HUR}_{\cdot,\cdot,\cdot}^m(D)$
(resp., $\text{HUR}_{\cdot,\cdot,\cdot}^{m,G}(D)$) of
$\text{HUR}_{d,b}^m(D)$ under the canonical map
$$\text{HUR}_{d,b}^m(D)\to
\text{HUR}_{d,b}(D)=\text{HUR}_{d,b}^m(D)/\mathcal S_d.$$ In
particular, we have $\text{HUR}_{d,s_1}(D)=\text{HUR}_{d,s_2}(D)$ if
and only if there is $\sigma\in\mathcal S_d$ such that
$\lambda(\sigma)(s_1)=s_2$.

Corollary \ref{uniq} gives us a complete description of irreducible
components of $\text{HUR}_{d,b}(D)$ in the case $d=3$.
\begin{cor} The
irreducible components of {\rm $\text{HUR}_{3,b}^G(D)$} are uniquely
defined by the monodromy factorization type and the type of global
monodromy if $G\simeq \mathcal S_2$ or $\mathcal S_3$.

The space {\rm $\text{HUR}_{3,b}^{A_3}(D)$} consists of
$m=[\frac{b}{6}]+1$ irreducible components if the global monodromy
is equal to $\bf{1}$ and it consists of $m=-[\frac{-b}{3}]$
irreducible components if the global monodromy is not equal to
$\bf{1}$.
\end{cor}
\subsection{Hurwitz spaces of (non-marked) coverings of $\mathbb P^1$}
In \cite{F}, Hurwitz spaces $\text{HUR}_{d,b}(\mathbb P^1)$ of
coverings of the projective line $\mathbb P^1$ of degree $d$,
branched over $b$ points, were described as non-ramified coverings
of the complement of the discriminant locus $\Delta$ in the
symmetric product $\mathbb P^{(b)}$ of $b$ copies of $\mathbb P^1$.
The choice of a point $\infty\in\mathbb P^1$ and the identification
$\mathbb C$ with $\mathbb P^1\setminus \{\infty\}$ defines an
imbedding of $\text{HUR}_{d,b}(D_{\infty})$ into
$\text{HUR}_{d,b}(\mathbb P^1)$ as an everywhere dense  open subset.
So we get the following
\begin{prop} \label{irred3} The irreducible components of {\rm
$\text{HUR}_{d,b}(\mathbb P^1)$} are in one to one correspondence
with the orbits of the action of $\mathcal S_d$ by simultaneous
conjugation on $\Sigma_{d,\bf{1}, b}=\{ s\in \Sigma_{d,\bf{1}}\mid
ln(s)=b\}$.
\end{prop}

As in subsection \ref{disc}, we can introduced the unions
$\text{HUR}_{\cdot,\cdot,\cdot}(\mathbb P^1)$ (resp.,
$\text{HUR}_{\cdot,\cdot,\cdot}^G(\mathbb P^1)$) of irreducible
components of $\text{HUR}_{d,b}(\mathbb P^1)$ for fixed elements of
$\Sigma_{b,\bf{1}}$, for fixed types of monodromy factorizations,
fixed Galois groups, and so on.

As a consequence of Proposition \ref{simple} we have
\begin{thm} There is a natural structure of the semigroup
$\Sigma_{d,{\bf{1}}}^{\mathcal S_d}=\{ s\in \Sigma_{d,{\bf{1}}} \mid
(\mathcal S_d)_s=\mathcal S_d\}$ on the set of irreducible
components of {\rm $\text{HUR}_{d}^{\mathcal S_d}(\mathbb P^1)$}.
\end{thm}

Theorem \ref{fuc1} and Corollary \ref{uniq} give us the following
two theorems.
\begin{thm} \label{TH1} The space {\rm $\text{HUR}_{d,t}^{\mathcal S_d}(\mathbb
P^1)$} is irreducible if the monodromy factorization type $t$
contains more than $3(d-1)$ transpositions.
\end{thm}

\begin{thm} \label{TH2} The irreducible components of {\rm $\text{HUR}_{3,b}^{G}(\mathbb P^1)$}
are uniquely defined by the monodromy factorization type if $G\simeq
\mathcal S_2$ or $\mathcal S_3$.

The space {\rm $\text{HUR}_{3,b}^{A_3}(\mathbb P^1)$} consists of
$m=[\frac{b}{6}]+1$ irreducible components.
\end{thm}

According to Theorems \ref{TH1}, \ref{TH2}, and Clebsch -- Hurwitz
Theorem, one can hope that the space $\text{HUR}_{d,t}^{\mathcal
S_d}(\mathbb P^1)$ is irreducible always for a fixed monodromy
factorization type $t$. The following theorem also confirms this
conjecture.

\begin{thm} \label{THM} Let $\sigma_1\in \mathcal S_d$ be a transposition
and $\sigma_2 \in \mathcal S_d$ be a cycle of length $d$. Then the
space {\rm $\text{HUR}_{d,t}(\mathbb P^1)$} is irreducible for fixed
type $t$ of the form $([2],t(\sigma_1\sigma_2^{-1}),[d])$.

There are exactly $[\frac{d}{2}]$ different types of such
form.\end{thm}

\proof If the type of $s\in \Sigma_d$ is
$([2],t(\sigma_2^{-1}\sigma_1),[d])$, then $ln(s)=3$ and hence
$\text{HUR}_{d,t}(\mathbb P^1)$ is unramified covering of $\mathbb
P^{(3)}\setminus \Delta$.

By Theorem \ref{TH2}, we can assume that $d\geq 4$.

Let us show that there are at least $[\frac{d}{2}]$ different
elements $s\in \Sigma_d$ of the form $s=x_{\sigma_1}\cdot
x_{\sigma_2}\cdot x_{\sigma_2^{-1}\sigma_1}$. For this it suffices
to show that there are $[\frac{d}{2}]$ different types for the
elements of $\mathcal S_d$ of the form $\sigma_2^{-1}\sigma_1$.
Indeed, without loss of generality, we can assume that
$\sigma_2^{-1}=(1,2)(2,3)\dots (d-1,d)$ and $\sigma_1= (i,d)$. Then
the type of
$$\begin{array}{ll}\sigma_2^{-1}\sigma_1= & (1,2)(2,3)\dots (d-1,d)(i,d)=
\\ & (1,2)\dots (d-2,d-1)(i,d-1)(d-1,d)=\dots =
\\ & (1,2)\dots (i-1,i)(i+1,i+2)\dots (d-1,d), \end{array}$$
is $[i,d-i]$ for $i=2,\dots ,[\frac{d}{2}]$ and $[d-1]$ for $i=1$.
In particular, the element $\sigma_2^{-1}\sigma_1$ is not conjugated
with $\sigma_1$ nor with $\sigma_2$ if $d\geq 4$.

Consider the set $U$ of words $w\in W$ consisting of three letters
$x_i$, $x_j$, $x_k$, where $x_i=x_{\sigma_1}$, $x_j=x_{\sigma_2}$,
and $x_k=x_{\eta}$, where $\eta$ is equal to  either
${\sigma_2^{-1}\sigma_1}$ or ${\sigma_1\sigma_2^{-1}}$ (depending on
the position of the letter $x_k$ in the word $w$ so to have
$\alpha(w)=\bf{1}$). Since the number of different transpositions is
equal to $\frac{d(d-1)}{2}$, the number of different cycles
$\sigma_2$ of length $d$ is equal to $(d-1)!$, and the element $x_k$
is uniquely defined by the positions of the letters $x_i$, $x_j$,
and $x_k$ in the word $w$ and by $\sigma_1$ and $\sigma_2$, then we
have
\begin{equation} \label{inq} \sharp U=6\frac{d(d-1)}{2}(d-1)!=3d!(d-1).\end{equation}

Consider two words $w_1$ and $w_2$ of $U$ consisting, respectively,
of letters $x_{i_1}=x_{\sigma_1}$, $x_{j_1}=x_{\sigma_2}$,
$x_{k_1}=x_{\eta}$ and $x_{i_2}=x_{\hat{\sigma}_1}$,
$x_{j_2}=x_{\hat{\sigma}_2}$, $x_{k_2}=x_{\hat{\eta}}$. It is easy
to see that the words $w_1$ and $w_2$ do not belong to the same
orbit of the action of $\mathcal S_d$ by simultaneous conjugation if
$t(\eta)\neq t(\hat{\eta})$. Therefore in $U$ there exist at least
$[\frac{d-1}{2}]$ different orbits of this action.  Let us fix a
word $w\in U$ and count the number of elements belonging to the
orbit of $w$. It is easy to see that the stabilizer of the letter
$x_{\sigma_2}$ is the cyclic subgroup $Z_{\sigma_2}$ of $\mathcal
S_d$ generated by $\sigma_2$. The transposition $\sigma_1$ is fixed
under the conjugation by $\sigma_2^n$ for $n\in [1, d-1]$ only if
$d=2n$ and in this case the order of the stabilizer of $w$ is less
or equal $2$. Like in the computation of the number of different
types of permutations of the form $\sigma_2^{-1}\sigma_1$, one can
show that if $d=2n$ and
$\sigma_2^{-n}\sigma_1\sigma_2^{n}=\sigma_1$, then $t(\eta)=[n,n]$.
We have
\begin{equation} \label{inq1} \sharp U\geq 6[\frac{d}{2}]d!=3d!(d-1)\end{equation}
if $d$ is odd and if $d=2n$ is even, then
\begin{equation} \label{inq2} \sharp U\geq 6([\frac{d}{2}]-1)d!+6\frac{d!}{2}=3d!((2n-1)=3d!(d-1).\end{equation}
It follows from (\ref{inq}) -- (\ref{inq2}) that the orbit under the
simultaneous conjugation of an element $s$ of type $\tau(s)=([2],
t(\sigma_2^{-1}\sigma_1), [d]\}$ is uniquely defined by its type.
Therefore  the space $\text{HUR}_{d,t}(\mathbb P^1)$ is irreducible
for fixed type $t=([2], t(\sigma_2^{-1}\sigma_1), [d])$ and the
number of such components is equal to $[\frac{d}{2}]$. \qed

\subsection{Hurwitz spaces of Galois coverings}
Let $f:C\to \mathbb P^1$ be a Galois covering with Galois group
$G=Gal(C/\mathbb P^1)$, that is, $G$ is the deck transformation
group of the covering $f$ and the quotient space $C/G=\mathbb P^1$.
In this case we have $\deg f=|G|$ and if we fix a point $\infty\in
\mathbb P^1$ over which $f$ is not ramified and fix a point $e\in
f^{-1}(\infty)$, then the action of $G$ on $f^{-1}(\infty)$ defines
a numbering of the points in $f^{-1}(\infty)$ by the elements of
$G$. If we choose also a numbering of the points in $f^{-1}(\infty)$
by the numbers belonging to the segment $I_{|G|}=[1,|G|]$, then
these numberings define an embedding $G\hookrightarrow \mathcal
S_{|G|}$. It is easy to see that this is Cayley's embedding.
Therefore the Hurwitz space $\text{HUR}^G(\mathbb P^1)$ of Galois
coverings with the Galois group $G$ can be identified with
$HUR_{|G|,{\bf{1}}}^G(\mathbb P^1)$ and, in particular, the natural
map
\begin{equation} \label{mmmmm} \text{HUR}_{|G|,{\bf{1}}}^{m,G}(\mathbb
P^1)\to \text{HUR}_{|G|,{\bf{1}}}^G(\mathbb
P^1)=\text{HUR}^G(\mathbb P^1)
\end{equation}
is surjective unramified morphism.
\begin{thm} \label{la} The irreducible components of {\rm $\text{HUR}^G(\mathbb P^1)$} are
in one to one correspondence with the orbits of the elements $s\in
S_G^G\subset S(G,G)$ under the action of $Aut(G)$ on $S(G,G)$.

If $Aut(G)=G$, then there is a natural structure of the semigroup
$S_{G,{\bf{1}}}^G$ on the set of irreducible components of {\rm
$\text{HUR}^G(\mathbb P^1)$}.
\end{thm}

\proof The first part of the theorem follows from  Corollary
\ref{norbits}.

To prove the second part, note that the equality $Aut(G)=G$ means
that any automorphism of $G$ is inner. By Proposition \ref{simple},
the elements of $S_{G,{\bf{1}}}^G$ are fixed under the action of $G$
by simultaneous conjugation. Therefore, by Corollary \ref{norbits},
natural map (\ref{mmmmm}) is an isomorphism which gives the desired
structure of semigroup on $\text{HUR}^G(\mathbb P^1)$. \qed

In particular, Theorem \ref{la} and Corollary \ref{uniq} imply

\begin{thm} The irreducible components of the Hurwitz space {\rm $\text{HUR}^{\mathcal S_3}(\mathbb P^1)$}
of Galois coverings with Galois group $G=\mathcal S_3$ are defined
uniquely by the monodromy factorization type of coverings belonging
to them.
\end{thm}

 \ifx\undefined\bysame
\newcommand{\bysame}{\leavevmode\hbox to3em{\hrulefill}\,}
\fi

\end{document}